\documentclass{jnmp}

\usepackage{amsmath}

\JNMPnumberwithin{equation}{section}

\begin{document}

\renewcommand{\evenhead}{F~Calogero and J-P Fran\c{c}oise}
\renewcommand{\oddhead}{Periodic Motions Galore}

\thispagestyle{empty}

\FirstPageHead{9}{1}{2002}{\pageref{calogero-firstpage}--\pageref{calogero-lastpage}}{Article}

\copyrightnote{2002}{F~Calogero and J-P Fran\c{c}oise}

\Name{Periodic Motions Galore: \\
How to Modify Nonlinear Evolution Equations\\
so that They Feature a Lot of Periodic Solutions}

\label{calogero-firstpage}

\Author{F~CALOGERO~$^\dag$ and J-P FRAN\c{C}OISE~$^\ddag$}

\Address{${}^{\dag}$~Dipartimento di Fisica, Universit\'a di Roma ``La Sapienza'',\\
$\phantom{{}^{\dag}}$~P. le A.~Moro 2, 00185 Roma, Italy,\\
$\phantom{{}^{\dag}}$~Istituto Nazionale di Fisica Nucleare, Sezione di Roma, Italy\\
$\phantom{{}^{\dag}}$~E-mail:
francesco.calogero@uniroma1.it, \ francesco.calogero@roma1.infn.it\\[10pt]
${}^{\ddag}$~GSIB, UFR de math\'ematiques, 175 rue du
Chevaleret,\\
$\phantom{{}^{\ddag}}$~Universit\'e de Paris VI, 75013 Paris, France\\
$\phantom{{}^{\ddag}}$~E-mail: jpf@ccr.jussieu.fr}

\Date{Received July 19, 2001; 
Accepted August 30, 2001}

\begin{abstract}
\noindent
A simple trick is illustrated, whereby nonlinear evolution equations
can be modified so that they feature a lot~-- or, in some cases,
only~-- {\it periodic} solutions. Several examples (ODEs and PDEs) are
exhibited.
\end{abstract}

\section{Introduction and outline of results}

Many natural phenomena (in physics, chemistry, biology, physiology, \ldots),
and also several phenomena pertaining to the social sciences (for
instance in economics), display a cyclic behavior. Hence, to the
extent they can be described by mathematical models, they call into
play evolution equations that feature {\it periodic} solutions.
Practitioners of such mathematical modeling (``applied mathematics'')
are therefore likely to profit from a simple trick, whereby evolution
equations can be modified in a rather neat way, so that they feature a
lot of~-- or, in some cases, only~-- {\it periodic} solutions.
Purpose and scope of this paper is to illustrate the efficacy of such
a trick, by applying it to several standard evolution equations, and
by tersely discussing the periodic behavior of solutions of these
modified equations. As we will see, this simple trick relates
periodicity in real time to analyticity in a (suitably introduced)
time variable, and thereby brings to light an explicit relation among
analyticity in complex time and integrable behavior of the
corresponding real dynamics. A representative list of such equations reads as
follows:

{\it First-order algebraic complex ODE}
\begin{equation}
\dot{w} - {\rm i} \Omega w = \alpha w^{p/q};
\end{equation}

{\it Polynomial vector field in the plane}
\begin{subequations}
\begin{gather}
\dot{u} + \Omega v = a_{1} U-a_{2} V, \qquad \dot{v} - \Omega u = a_{1}
V + a_{2} U , \\
U\equiv \sum_{m=0}^{[[p/2]]}(-1)^{m} \left(\begin{array}{@{}c@{}} p\\ 2m\end{array}\right) u^{p-2m}
v^{2m}, \nonumber\\ V \equiv \sum_{m=0}^{[[(p-1)/2]]} (-1)^{m}
\left(\begin{array}{@{}c@{}}  p\\ 2m+1\end{array}\right) u^{p-2m-1} \, v^{2m+1}; \\
\dot{u} + \Omega v = a_{1} \left(u^{2} - v^{2}\right) -2 a_{2} uv, \qquad \dot{v} -
\Omega u = a_{2} \left(u^{2}-v^{2}\right) + 2 a_{1} uv; \tag{1.3}
\end{gather}
\end{subequations}

{\it Oscillator with additional cubic force}
\setcounter{equation}{3}
\begin{equation}
\ddot{x} + (\Omega /2)^{2} x = a^{2} x^{-3};
\end{equation}

{\it Painlev\'e-type equations (complex and real versions)}
\begin{gather}
\ddot{w} + \Omega^{2} w= \left(\alpha w^{2} + \gamma\right) \exp (5 {\rm i} \Omega
t); \\
\ddot{w} - 5{\rm i} \Omega \dot{w} - 6 \Omega^{2} w = \alpha w^{2} +
\gamma \exp (5 {\rm i} \Omega t); \\
\ddot{w}  +5{\rm i} \Omega \dot{w} - 6 \Omega^{2} w = \alpha w^{2} \exp
 (5{\rm i} \Omega t) + \gamma; \\
\ddot{w} -3 {\rm i} \Omega \dot{w} - 2 \Omega^{2}w= \alpha w^{3}+ (\gamma
w + \delta ) \exp (3 {\rm i} \Omega t); \\
\ddot{u} + \Omega^{2} u = \cos (5 \Omega t) \left[a_{1} \left(u^{2} - v^{2}\right) -
2 a_{2} uv +c_{1}\right]\nonumber\\
\qquad{}- \sin (5 \Omega t) \left[a_{2} \left(u^{2} - v^{2}\right) + 2
a_{1}uv + c_{2}\right], \tag{1.9a}
\\
\ddot{v} + \Omega^{2} v = \sin (5 \Omega t) \left[ a_{1} \left(u^{2}- v^{2}\right) -
2a_{2} u v + c_{1} \right]\nonumber\\
\qquad{}+ \cos (5 \Omega t)
\left[ a_{2} \left(u^{2} - v^{2}\right) + 2 a_{1}u v + c_{2}\right]; \tag{1.9b}
\\
\ddot{u} +5 \Omega \dot{v} - 6 \Omega^{2}u = a_{1} \left(u^{2}- v^{2}\right) -
2a_{2} u v + c_{1} \cos (5 \Omega t) - c_{2}\sin (5 \Omega t),
\tag{1.10a}\\
\ddot{v} -5\Omega \dot{u} - 6 \Omega^{2} v = a_{2}\left(u^{2}- v^{2}\right) +
2 a_{1} u v + c_{2} \cos (5 \Omega t)  + c_{1} \sin (5 \Omega
t); \tag{1.10b}\\
\ddot{u} - 5 \Omega \dot{v} - 6 \Omega^{2} u = \left[a_{1}\left(u^{2}-v^{2}\right)
-2a_{2} uv\right] \cos (5 \Omega t) \nonumber\\
\qquad{}- \left[a_{2} \left(u^{2}- v^{2}\right) + 2 a_{1} u
v\right] \sin (5 \Omega t)+c_{1}, \tag{1.11a}\\
\ddot{v}+5 \Omega \dot{u} - 6 \Omega^{2}v = \left[a_{1} \left(u^{2} -v^{2}\right)
-2a_{2}uv\right]\sin (5 \Omega t) \nonumber\\
\qquad{}+ \left[a_{2} (u^{2}- v^{2}) + 2a_{1} uv\right]
\cos (5\Omega t) +c_{2}; \tag{1.11b}\\
\ddot{u} + 3 \Omega \dot{v} - 2\Omega^2 u = a_1 u\left(u^2 - 3 v^{2}\right) -
a_{2} v\left(3 u^{2} - v^{2}\right)\nonumber\\
\qquad{}
+ \left(c_{1} u - c_{2} v + d_{1}\right) \cos (3 \Omega t) - (c_{2} u + c_{1}
v + d_{2}) \sin (3 \Omega t), \tag{1.12a}\\
\ddot{v} - 3 \Omega \dot{u} - 2 \Omega^{2} v = a_{2} u \left(u^{2} - 3
v^{2}\right) + a_{1} v\left(3u^{2} - v^{2}\right)\nonumber\\
\qquad{}
+ \left(c_{1} u - c_{2} v + d_{1}\right) \sin (3\Omega t) + (c_{2} u + c_{1}
v  + d_{2}) \cos (3\Omega t);\tag{1.12b}
\end{gather}

{\it Autonomous second-order ODEs (complex and real versions)}
\setcounter{equation}{12}
\begin{gather}
\ddot{w} - {\rm i} (5/2) \Omega \dot{w} - (3/2) \Omega^{2} w = \alpha
w^{2};\\ \ddot{w} - {\rm i} \Omega \dot{w} - 2 \Omega^{2} = \alpha
\exp (w);\tag{1.14a}\\ \ddot{w}-{\rm i}\Omega \dot{w}  + 2
\Omega^{2} w =(\dot{w} - {\rm i} \Omega w) w ;\tag{1.14b}\\
\ddot{w} - 3 {\rm i} \Omega \dot{w} -2 \Omega^{2} w = (\dot{w} - {\rm i} \Omega
w) w;\tag{1.15}\\
 \ddot{u} + (5/2) \Omega \dot{v} - (3/2)
\Omega^{2} u = a_{1} \left(u^{2} - v^{2}\right) - 2a_{2}
uv,\tag{1.16a}
\end{gather}
\begin{gather}
 \ddot{v} - (5/2) \Omega \dot{u} - (3/2) \Omega^{2}
v = a_{2} \left(u^{2} - v^{2}\right) + 2a_{1} uv;\tag{1.16b}\\
\ddot{u} + \Omega \dot{v} - 2 \Omega^{2}=  \exp (u) [a_{1} \cos
(v) - a_{2} \sin (v)], \tag{1.17a}\\ \ddot{v} - \Omega \dot{u} =
\exp (u) [a_{2}\cos (v) + a_{1}  \sin (v)];\tag{1.17b}\\ \ddot{u}
+ \Omega \dot{v} +2 \Omega^{2}u= \dot{u} u - \dot{v} v + 2 \Omega
uv, \tag{1.18a}\\ \ddot{v} - \Omega \dot{u} + 2 \Omega^{2}v =
\dot{u} v +u \dot{v} - \Omega \left(u^{2}-
v^{2}\right);\tag{1.18b}\\ 
\ddot{u} + 3 \Omega \dot{v} - 2
\Omega^{2} u = \dot{u} u - \dot{v} v + 2 \Omega uv,\tag{1.19a}\\
\ddot{v} - 3 \Omega \dot{u} - 2 \Omega^{2} u = \dot{u} v + u
\dot{v} - \Omega \left(u^{2} - v^{2}\right);\tag{1.19b}\\
\ddot{w}+ {\rm i} \{(3q_{1}q_{2}+ p_{1}q_{2}-p_{2}q_{1}) /[q_{1}
(p_{2}-2q_{2})]\} \Omega \dot{w}\nonumber\\ \qquad{} + [(q_{2}/
q_{1}) (p_{1} + q_{1}) / (p_{2} - 2q_{2})] \Omega^{2} w = \alpha
(\dot{w} - {\rm i} \Omega w)^{p_{2}/q_{2}} w^{p_{1}/q_{1}}
;\tag{1.20}\\ \ddot{w} - 3 {\rm i}\Omega \dot{w} - 2 \Omega^{2}w =
\alpha (\dot{w} - {\rm i}\Omega w)^{3} w^{-3}; \tag{1.21}\\
\ddot{w} + {\rm i} (nm + n-1) \Omega \dot{w} + n(m+1)\Omega^{2}w =
\alpha (\dot{w} - {\rm i} \Omega w)^{(2n+1)/n} w^{m};\tag{1.22}\\
\ddot{w}- {\rm i}~ (2nm +1) \Omega \dot{w}- 2nm \Omega^{2}w
 = \alpha (\dot{w} - {\rm i} \Omega w)^{(2n+1)/n}
w^{-(2m+1)};\tag{1.23}\\
\ddot{w}-3 {\rm i}\Omega \dot{w} -2 \Omega^{2} w= \alpha
w(\dot{w}-{\rm i}\Omega w) + \beta w^{3};\tag{1.24}\\
\ddot{w} -3 {\rm i} \Omega \dot{w}- 2 \Omega w= \alpha w (\dot{w} -
{\rm i} \Omega
w) - (\alpha /3)^{2} w^{3}; \tag{1.25}\\
\ddot{w} -3 {\rm i} \Omega \dot{w} -2 \Omega^{2} w= \alpha w
(\dot{w} - {\rm i}
\Omega w);\tag{1.26}\\
\ddot{u} +3 \Omega \dot{v} - 2 \Omega^{2}u = a_{1} (u \dot{u} - v
\dot{v} + 2 \Omega uv) \nonumber\\
\qquad{}- a_{2}\left[\dot{u}v + u \dot{v} - \Omega
\left(u^{2}-v^{2}\right)\right]
 + b_{1} u \left(u^{2}-3v^{2}\right) - b_{2}v \left(3u^{2}-v^{2}\right), \tag{1.27a}\\
\ddot{v} - 3 \Omega \dot{u} - 2\Omega^{2} v = a_{2}(u \dot{u} - v
\dot{v} + 2 \Omega uv)\nonumber\\
\qquad{} + a_{1}\left[\dot{u} v+ u \dot{v} - \Omega \left(u^{2} -
v^{2}\right)\right]+ b_{2} u \left(u^{2}- 3v^{2}\right) + b_{1}v\left(3 u^{2}-v^{2}\right);
\tag{1.27b}\\
\ddot{w} - {\rm i} \Omega \dot{w} = \dot{w}^{2} F(w); \tag{1.28}
\end{gather}

{\it Autonomous third order ODEs (complex and real versions)}
\setcounter{equation}{28}
\begin{gather}
\dddot{w}-10{\rm i}\Omega \ddot{w} - 31 \Omega^{2}\dot{w}+30{\rm
i} \Omega^{3}w = \alpha (2 \dot{w} -5 {\rm i}\Omega
w)w;\\
\dddot{z}-10{\rm i} \Omega \ddot{z} - 19 \Omega^{2} \dot{z} -
30 {\rm i} \Omega^{3}z=\alpha (2 \dot{z} - 5 {\rm i} \Omega z)z;
\\
\dddot{w}+31 \Omega^{2} \dot{w}+ 30 {\rm i}\Omega^{3}w+ 5{\rm
i}\Omega \gamma = 2(\ddot{w} + 5 {\rm i} \Omega \dot{w} - \gamma)
\dot{w}/w; \\
\dddot{w} -5{\rm i} \Omega \ddot{w} + \Omega^{2}\dot{w} -5 {\rm
i} \Omega^{3}w = 2 w \dot{w} \left(\ddot{w}+ \Omega^{2}w\right) / \left(w^{2}+ \eta\right);
\\
\dddot{u}+10\Omega \ddot{v} -31 \Omega^{2} \dot{u} -
30 \Omega^{3} v\nonumber\\
\qquad {}=2a_{1}(\dot{u} u - \dot{v}v + 5 \Omega u v) -a_{2}\left[2(\dot{u} v+ \dot{v}
u)-5 \Omega \left(u^{2}-v^{2}\right)\right],\tag{1.33a}\\
\dddot{v}-10\Omega \ddot{u} -31 \Omega^{2} \dot{v} +
30 \Omega^{3} u\nonumber\\
\qquad{}=2a_{2}(\dot{u} u - \dot{v}v + 5 \Omega u v) + a_{1}\left[2(\dot{u} v+ \dot{v}
u)-5 \Omega \left(u^{2}-v^{2}\right)\right];\tag{1.33b}\\
\dddot{u}+10\Omega \ddot{v} -19 \Omega^{2} \dot{u} +
30 \Omega^{3} v\nonumber\\
\qquad{}=2a_{1}(\dot{u} u - \dot{v}v + 5 \Omega u v) -a_{2}\left[2(\dot{u} v+ \dot{v}
u)-5 \Omega \left(u^{2}-v^{2}\right)\right],\tag{1.34a}\\
\dddot{v}-10\Omega \ddot{u} -19 \Omega^{2} \dot{v} -
30 \Omega^{3} u\nonumber\\
\qquad {}=2a_{2}(\dot{u} u - \dot{v}v + 5 \Omega u v) -a_{1}\left[2(\dot{u} v+ \dot{v}
u)-5 \Omega \left(u^{2}-v^{2}\right)\right];\tag{1.34b}
\end{gather}

{\it First-order PDE of shock type  (complex and real versions)}
\setcounter{equation}{34}
\begin{gather}
w_{t}- {\rm i} \Omega w= \alpha w_{x} w^{p/q};
\end{gather} \begin{gather}
u_{t}+\Omega v = a_{1}(u_{x}u - v_{x}v)
-a_{2}(u_{x}v+v_{x}u),\tag{1.36a} \\
v_{t}-\Omega u = a_{2}(u_{x}u
- v_{x}v) +a_{1}(u_{x}v+v_{x}u);\tag{1.36b}\\ u_{t}+\Omega v
=a_{1}\left[u_{x}\left(u^{2}-v^{2}\right) -2v_{x}uv\right]
-a_{2}\left[2u_{x}uv+ v_{x}\left(u^{2}-v^{2}\right)\right],
\tag{1.37a}\\ v_{t}-\Omega v
=a_{2}\left[u_{x}\left(u^{2}-v^{2}\right) -2v_{x}uv\right]
+a_{1}\left[2u_{x}uv+ v_{x}\left(u^{2}-v^{2}\right)\right];
\tag{1.37b}
\end{gather}

{\it Burgers-type equations  (complex and real versions)}
\setcounter{equation}{37}
\begin{gather}
w_{t}= \left({\rm i} \Omega x w + \beta w_{x} + \alpha w^{2}\right)_{x}; \\
u_{t}=\left[-\Omega x v + b_{1}u_{x}-b_{2}v_{x}+ a_{1}\left(u^{2}-v^{2}\right) - 2
a_{2}uv\right]_{x}, \tag{1.39a}\\
v_{t}=\left[\Omega x u + b_{2}u_{x}-b_{1}v_{x}+ a_{1}\left(u^{2}-v^{2}\right) - 2
a_{1}uv\right]_{x}; \tag{1.39b}
\end{gather}

{\it Generalized KdV-type equations}
\setcounter{equation}{39}
\begin{gather}
w_{t} = {\rm i} \Omega [2(q/p) w + x w_{x}] + \beta w_{xxx}+ \alpha w_{x}
w^{p/q};
\end{gather}

{\it KdV-type equation (complex and real versions)}
\begin{gather}
w_{t} = {\rm i} \Omega [2 w + x w_{x}] + \beta w_{xxx}+ \alpha w_{x} w;\\
u_{t}=-\Omega (2v +xv_{x}) +b_{1} u_{xxx}-b_{2} v_{xxx}+ a_{1}
(u_{x}u - v_{x}v) - a_{2}(u_{x}v+ v_{x}u), \tag{1.42a}\\
v_{t}= \Omega (2 u+ x u_{x}) + b_{1} v_{xxx}+ b_{2}  u_{xxx}+
a_{2}(u_{x}u - v_{x}v) + a_{1}(u_{x}v+v_{x}u); \tag{1.42b}
\end{gather}

{\it Modified KdV-type equation (complex and real version)}
\setcounter{equation}{42}
\begin{gather}
w_{t}= [{\rm i} \Omega x w +\beta w_{xx}+ (\alpha/3) w^{3}]_{x}; \\
u_{t}=\left[-\Omega x v + b_{1} u_{xx} - b_{2} v_{xx} +(a_{1}/3) u
\left(u^{2}-3 v^{2}\right) - (a_{2}/3) v \left(3 u^{2}-v^{2}\right)\right]_{x}, \tag{1.44a}\\
v_{t}=\left[\Omega x u + b_{2} u_{xx} + b_{1} v_{xx} +(a_{2}/3) u
\left(u^{2}-3 v^{2}\right) + (a_{1}/3) v \left(3 u^{2}-v^{2}\right)\right]_{x}; \tag{1.44b}
\end{gather}

{\it KP-type equation (complex and real versions)}
\setcounter{equation}{44}
\begin{gather}
[w_{t} - {\rm i} \Omega w - ({\rm i}/2) \Omega x w_{x}-{\rm i}
\Omega  y w_{y}+ \beta w_{xxx} +
\alpha w_{x}  w]_{x} + \gamma w_{yy} = 0;\\
[u_{t} + \Omega v + (\Omega/2) x v_{x} + \Omega y v_{y} + b_{1}
u_{xxx} - b_{2} v_{xxx}\nonumber\\
\qquad {}+ a_{1} (u_{x} u - v_{x} v) -a_{2} (u_{x} v + v_{x} v)]_{x} + c_{1}
u_{xx} -c_{2} v_{yy} = 0,\tag{1.46a}\\
[v_{t} - \Omega u - (\Omega/2) x u_{x} - \Omega y u_{y} + b_{1}
v_{xxx} + b_{2} x_{xxx} + b_{2} u_{xxx}\nonumber\\
\qquad {}+ a_{2} (u_{x} u - v_{x} v) + a_{1} (u_{x} v + v_{x} u)]_{x} + c_{2}
u_{yy} + c_{1} v_{yy} = 0.\tag{1.46b}
\end{gather}

These evolution equations are discussed tersely, one by one, at the
end of this introductory Section~1. But let us immediately emphasize
that the independent variable,~$t$ (``time''), which accounts for their
evolutionary character, is of course assumed to be {\it real}, and
also {\it real} is the nonvanishing constant $\Omega$, to which the
period
\setcounter{equation}{46}
\begin{gather}
T= 2 \pi/|\Omega|
\end{gather}
is associated. The dependent variables are in general complex, as are
other constants featured by these equations; but sometimes the
evolution equations are written in the real form  yielded by an
explicit separation of all complex numbers into their real
and imaginary parts. The main point is
that {\it all} these evolution equations (but one, namely (1.14a) and
its real version (1.17); see below) feature {\it many}, or in
some cases {\it only}, solutions which are periodic with period $T$~--
or, in several cases, with periods obtained by multiplying $T$ by a
rational number. Note that the above list includes both ODEs and PDEs;
in the latter case we restricted (with one exception, see (1.45), 
(1.46)),
merely for the sake of simplicity,
consideration to (few) equations in $2 = 1+1$ variables, namely to only one,
of course {\it real}, ``space'' variable, $x$, in addition to the,
also {\it real}, ``time'' variable~$t$. The alert reader, after
having understood how the trick which yielded these evolution
equations works, will have no difficulty (rather, some fun!) in
applying it more widely.

The trick is introduced in the following Section~2 to derive and
discuss the three prototypical evolution equations (1.1), (1.35) and
(1.38); it is then used in Section 3 to derive and discuss the other
evolution equations reported in the above list (as well as some more
general versions encompassing several of them).

The elementary character of the trick presumably precludes any hope to ascertain
{\it precisely} who introduced it firstly~[1,~2]. An ample use of it, in the
context of classical many-body problems amenable to exact treatments,
has been made in Ref.~[3], to which we refer for references to
previous applications in such a context; for this reason we did not
include in the above list of evolution equations any (system of) ODEs
describing classical many-body motions.

Let us end this introductory Section~1 by providing some details on the
evolution equations listed above.

In the ODE (1.1), as well as in all subsequent ODEs in the list, dots
denote differentiations with respect to the ({\it real}) independent
variable $t$ (``time''). In this ODE~(1.1), the dependent variable
$w \equiv w (t)$ is {\it complex}, $\alpha$ is a nonvanishing, but
otherwise arbitrary (possibly {\it complex}) constant, and $p/q$ is a
{\it rational} number different from unity, $p\neq q$. {\it All}
{\it nonsingular} solutions of this ODE, (1.1), are {\it periodic}: for a
given assignment of the rational number $p/q$, there are generally
two periods, $T_{1}$ and $T_{2}$, both integer multiples of~$T$, see
(1.47), and the initial data, $w(0)$, are divided into two
complementary open sets, which yield solutions periodic with one,
respectively with the other, period; these two sets of initial data
are separated by (a lower dimensional set of) initial data that yield
singular solutions. The solutions of this ODE, (1.1), are given, and
discussed, in Section~2.

In the system of two coupled ODEs (1.2) the two dependent variables,
$u \equiv u (t) $, $v \equiv v (t)$, are instead {\it real}, as well as
the two, otherwise {\it arbitrary}, constants $a_{1}$, $a_{2}$; $p$ is
an {\it integer} larger than unity, $p >1$; and the symbol $[[\cdots]]$
denotes the integral part of the doubly-bracketed quantity. This
system, (1.2), is merely another avatar of (1.1) with $q=1$, obtained
by setting $w = u + {\rm i} v$, $\alpha = a_{1} + {\rm i} a_{2}$.

The system (1.3) is just (1.2) with $p=2$.

The ODE  (1.4) is an avatar of (1.3) with $a_{1}=0$, $a_{2}=a$,
obtained by time-differentiating the first of the (1.3), by using the
second of the (1.3) to eliminate $\dot{v}$, by then using the first
of the (1.3) to eliminate $v$, and finally by setting $u(t) = [x
(t)]^{-2} - \Omega /(2a)$; the fact that all (real) solutions of this real
ODE are periodic with period $T$, see (1.47), is of course a
well-known result.

In the 3 complex evolution ODEs (1.5), (1.6) and (1.7) the two
constants $\alpha$ and $\gamma$ are {\it arbitrary} (possibly
complex), and also {\it arbitrary} (possibly complex) are the three
constants $\alpha, \gamma$ and $\delta$, in the complex ODE (1.8).
{\it All (nonsingular)} solutions of
these 4 (nonautonomous!) evolution ODEs are periodic with period $T$,
see (1.47).

The 4 systems of 2 {\it real} coupled evolution ODEs (1.9), (1.10),
(1.11) respectively (1.12) are merely the {\it real} avatars of
(1.5), (1.6), (1.7) respectively (1.8), obtained by setting $w = u +
{\rm i} v$, $ \alpha = a_{1} + {\rm i}a_{2}$,
$\gamma = c_{1} + {\rm i} c_{2}$, $\delta = d_{1} +
{\rm i}d_{2}.$ Of course {\it
all} their ({\it nonsingular}) solutions are periodic with period
$T$, see (1.47).

{\it All nonsingular} solutions of (1.13), (1.14b) and (1.15) are periodic in $t$ with
period (at most) $T$, see (1.47). The general solutions of the evolution
ODEs (1.13) and (1.14b) are displayed in Section~3, as well as conditions
sufficient to guarantee their {\it nonsingularity} (see (3.17)).
The solutions of (1.14b) coincide with the {\it time-derivative} of the
solutions of (1.14a).

The four systems of 2 {\it real} coupled evolution ODEs (1.16), (1.17), (1.18), (1.19)
are merely the {\it real} avatars of (1.13), (1.14a), (1.14b), (1.15), obtained by setting
$w = u + {\rm i}v$, $ \alpha = a_{1} +
{\rm i}a_{2}$; hence their ({\it real}) solutions can be immediately
obtained from the solutions of the complex ODEs (1.13), (1.14a),
(1.14b), (1.15), see immediately above.

The similarity among (1.14b) and (1.15) (and likewise among (1.18) and
(1.19)) should be noted; also note that (1.15) is a special case of
(1.20), which is not the case of (1.14b).

The complex evolution equation (1.20), with $\alpha$ an {\it
arbitrary} (possibly complex) constant and $p_{1}$, $q_{1}$, $p_{2}$,
$q_{2}$ four {\it integers} (but only the two rational numbers $p_{1} /
q_{1}$, $p_{2} / q_{2}$ actually enter in (1.20)), is expected to
possess lots of periodic solutions, as entailed by its derivation in
Section~3. However, only the 3~special cases of (1.20) corresponding
to (1.21), (1.22) and (1.23) are treated in any detail in Section~3
(in addition to (1.15), see above).
In these 3 evolution ODEs, (1.21)--(1.23), $\alpha$ is again an arbitrary,
possibly complex, constant; in (1.22) and (1.23) $n$ is a {\it
positive} integer, while $m$ is a {\it nonnegative} integer in
(1.22), a {\it positive} integer in (1.23). In all three cases,
(1.21)--(1.23), {\it all nonsingular} solutions are periodic. In the case
of (1.21), the general solution is exhibited in Section~3, and it is
shown there that in this case the {\it nonsingular} solutions split into two
sets, both however periodic with the {\it same} period~$T$,  see
(1.47); and these two sets are separated by a lower
dimensional set of initial data (characterized by (3.23) with (3.22)),
to which there correspond {\it singular} solutions of (1.21), namely
solutions such that $\dot{w} (t)$ diverges at {\it real} times
$t=t_{b}$ (with $t_{b}$ defined mod($T$) by (3.24) with (3.22) and
(3.23)).

In the other two cases, (1.22), (1.23), the situation is analogous, but
richer. Depending on the values of the integers $n$ and $m$, many
more periodicities are possible: and again these different
periodicities correspond to different sets of initial data, $w (0)$,
$\dot{w} (0)$, separated by lower-dimensional sets of such data which
themselves yield singular solutions (for which $\dot{w} (t)$ diverges
at some {\it real} times $t = t_{b}$, generally defined mod$\,(t_{p})$,
see (2.2)). More details in Section~3, where expressions of the
general solution of these two evolution ODEs, (1.22), (1.23), are
provided, albeit in somewhat implicit form.

The complex evolution ODE (1.24) features the two {\it arbitrary}
(possibly complex) constants $\alpha$, $\beta$, in addition  to the
{\it real} constant $\Omega$; its derivation in Section~3 entails that
it possesses lots of periodic solutions. Indeed the next two ODEs
listed, (1.25) respectively (1.26), which are clearly two special
cases of (1.24) (corresponding to $\beta = - (\alpha /3)^{2}$
respectively $\beta=0$), are explicitly solved in Section~3 (see (3.38) respectively
(3.41)), and it is shown there that {\it all} their {\it
nonsingular} solutions are {\it periodic} in $t$ with period~$T$, see
(1.47) (conditions necessary and sufficient to guarantee that these
solutions be {\it nonsingular} for all {\it real} values of $t$ are
also provided there, see (3.39) respectively (3.42)).

The system of two {\it real} evolution ODEs (1.27) is merely the real
avatar of (1.24), obtained by setting $w = u + {\rm i} v$, $\alpha = a_{1} +
{\rm i} a_{2}$, $\beta = b_{1} + {\rm i}b_{2}$.

In the second-order evolution ODE (1.28) $F(w)$ is a largely
arbitrary (but analytic) function; the fact that this evolution ODE,
in spite of its generality (associated with the large arbitrariness
of the function $F(w)$) possesses lots of periodic solutions is,
however, not new~[2], and we therefore forsake
any further discussion of this ODE in this paper.

The {\it third-order} complex ODEs (1.29) respectively (1.30) are merely
two avatars of the second order ODE (1.6), and the systems of two
real ODEs (1.33) respectively (1.34) are merely the real versions of
(1.29) respectively (1.30), obtained by setting $w(t)=u(t) +{\rm
i}v(t)$ respectively $z(t)=u(t)+{\rm i}v(t)$ and $\alpha=a_{1}+{\rm
i}a_{2}$; likewise the {\it third-order} ODEs (1.31) respectively
(1.32) are avatars of (1.7) respectively (1.5).
These facts are proven  in Section~3; they entail that {\it all
nonsingular} solutions of these ({\it third-order, autonomous}) ODEs,
(1.29)--(1.34), are {\it periodic} with period $T$, see (1.47).

In the PDE (1.35), as well as in all subsequent PDEs in the list, the
independent variables, $x$ and $t$, are supposed to be {\it real};
the independent variable, $w \equiv w (x,t)$, is generally {\it complex}
(but we also exhibit in some cases the {\it real} avatars of these
evolution PDEs: for instance (1.36) respectively (1.37) are the {\it
real} avatars of (1.35) with $q=1$ and $p=1$ respectively $p=2$, obtained by
setting $w = u+{\rm i}v$, $\alpha = a_{1} + {\rm i}a_{2}$). Of course subscripted
variables, here and below, denote partial differentiations. The
constant $\alpha$ in (1.35) is {\it arbitrary}, possibly complex;
$\Omega$ is (as always) a {\it real} (nonvanishing!) constant, and~$q$, $p$ are two
arbitrary {\it integers}. The initial-value problem for this
evolution PDE, (1.35), is solved in Section~2, and it is proven there
that there exist classes of initial data, $w(x,0)= w_{0} (x)$, which
yield solutions periodic in $t$ (for all values of $x$).

Likewise, it is shown in Section~2 that the PDE (1.38), with $\alpha$
and $\beta$ two {\it arbitrary} (possibly complex) constants and of
course $\Omega$ a {\it real} (nonvanishing!) constant, also possesses
(nonsingular!) solutions which are periodic in the ({\it real})
independent variable $t$ for all values of the ({\it real})
independent variable~$x$. Note the explicit appearance, in this PDE,
(1.38), of the variable~    $x$, which entails that this evolution PDE,
(1.38), is {\it not} translation-invariant. Clearly (1.39) is the
{\it real} avatar of (1.38), obtained by setting $w =u +{\rm i}v$, $\alpha =
a_{1} + {\rm i} a_{2}$, $\beta = b_{1} + {\rm i}b_{2}$.

The remaining evolution PDEs of the list reported above have been
derived using the trick, as described in detail in Section~3;
or they are special cases, or real avatars, of previous equations
in the list. It is therefore justified to expect that they all possess
lots of periodic solutions; but this hunch is not substantiated, in
Section~3, via the explicit exhibition of solutions, nor is it
discussed in any more detail here (lest the length of this paper
become excessive). Hence a more thorough investigation of these
evolution PDEs, as well as of the many others that can be easily
manufactured using the trick described in detail in Section~2, remains
as a task for the future.

Let us end this introductory Section~1 by indicating that our
presentation, as it unfolded above and as it continues below, has been
deliberately adjusted to cater not only to the interest of applied
mathematicians, but as well to the needs of scientists whose
mathematical training is somewhat superficial, because we believe
there are such potential customers who might usefully take advantage,
in various applicative contexts, of the findings reported in this paper.
This may have entailed occasional repetitions (for instance in the guise
of presenting both the complex
and the real versions of certain evolution equations), as well as the
rather systematic practice of keeping, when writing evolution
equations, also constants which might be
eliminated by simple rescalings of (dependent and independent)
variables. We hope most readers will agree that the advantages of such
an approach outweigh its defects.

\section{The trick}

The idea is to start from an evolution equation, chosen as candidate
for a modification the goal of which is to obtain a (modified)
evolution equation that features (a lot of) periodic solutions; to
then assume that the {\it independent} variable, $\tau$, in terms of
which the (original, unmodified) evolution unfolds, is {\it complex};
and to then introduce the {\it real} ``time'' variable,~$t$, so
that $\tau = \tau (t)$ be a {\it periodic} function of $t$, say
\begin{gather}
\tau \equiv \tau (t) = [\exp ({\rm i} \omega t) -1]/({\rm i} \omega). \tag{2.1a}
\end{gather}
Hereafter $\omega$ is a {\it positive} number, and we set
\setcounter{equation}{1}
\begin{gather}
t_{p} = 2 \pi /\omega,
\end{gather}
so that $\tau (t)$ is periodic in $t$ with period $t_{p}$. Indeed as $t$
varies over one period, say from $t=0$ to $t=t_{p}$, $\tau$ travels
from $\tau=0$ back to $\tau=0$, over a circular contour, $C$, in the
complex $\tau$-plane, centered at $\tau={\rm i}/ \omega$ and of radius
$r=1/ \omega= t_{p}/(2 \pi)$ (so that the diameter of this circle, $C$,
has length $2/\omega$ and straddles the upper imaginary axis in the
complex $\tau$-plane, from the origin, $\tau =0$, to $\tau = 2{\rm i}/
\omega$). Hence any function of the complex variable $\tau$, say
$\varphi (\tau)$, which is holomorphic in $\tau$ in a region of the
complex $\tau$-plane that includes the circle $C$, is a periodic
function of $t$ with period $t_{p}$. But if $\varphi (\tau)$ satisfies an
evolution equation, say
\begin{gather}
\varphi'  ( \tau) = F[\varphi (\tau), \tau],\qquad
\varphi'  ( \tau) \equiv d \varphi ( \tau) /d \tau,
\end{gather}
with $F(\varphi , \tau)$ an {\it analytic} function of its two
arguments, any solution $ \varphi ( \tau)$ of (2.3) characterized by an
initial datum $\varphi (0)$ such that $F(\varphi, \tau)$ is
holomorphic in the neighborhood of $\varphi = \varphi (0)$,
$\tau=0$, is itself guaranteed (by the standard
existence/uniqueness/analyticity theorem for evolution equations)
to be a holomorphic function of
$\tau$ in a disk $D$ centered, in the complex $\tau$-plane, at
$\tau=0$, and having a nonvanishing radius $\rho$ whose magnitude
depends on the function $F(\varphi, \tau)$ and on the initial datum
$\varphi(0)$. Clearly if the radius $\rho$ exceeds the diameter
$2/\omega$ of the circle $C$, $\rho > 2/\omega$, the disk $D$
includes the circle $C$, hence the corresponding solution $\varphi
(\tau)$, considered as a function of the {\it real} variable $t$, is
then {\it periodic} with period $t_{p}$, see (2.2). And it is of
course easy, via (2.1), to recast (2.3) as an evolution equation in
terms of the ({\it real}) time $t$ rather than the ({\it complex})
independent variable $\tau$.

At this stage one might wonder why was it necessary to let $\tau
\equiv\tau (t)$ be {\it complex}, see (2.1a); could not one instead
choose $\tau$ to be a {\it real} periodic function of $t$, say $\tau =
\omega^{-1} \sin (\omega t)$ or $\tau = \omega^{-1} {\rm tan} (\omega t)$? Of
course one could. But the additional, important requirement for the
basic idea to be {\it usefully} applicable is that the modified
evolution equation obtained via this approach in terms of the ({\it
real}) independent variable $t$ (``time'') have a {\it neat} look.
This is essential for such an equation to be eventually identified as
the appropriate evolution equation to model some natural phenomenon.
To this end the change of variable (2.1a) appears particularly
suitable, as the examples exhibited in this paper show; especially if
it is associated with an appropriate change of {\it dependent}
variable, typically by setting
\begin{gather}
w (t) = \exp ({\rm i} \lambda \omega t) \varphi [\tau (t)],
\end{gather}
of course with $\tau (t)$ defined by (2.1a) and with $\lambda$ a ({\it
rational}) number to be chosen conveniently (see below). And a
similar kind of modification can be conveniently applied to other
(``space'') variables in the case of PDEs, although this may also
cause certain difficulties, as we shall see below.

Enough now of this introductory discussion, which was merely meant to
convey the gist of this approach. Let us rather proceed and show how
this technique works, by demonstrating its effectiveness in three
(very simple, completely solvable) examples, one ODE and two PDEs.
But before doing so, let us note that (2.1a) entails
\begin{gather}
\dot{\tau} (t) = \exp ({\rm i} \omega t),\tag{2.1b}
\end{gather}
as well as
\begin{gather}
\tau (0) =0,\tag{2.1c}\\
\dot{\tau} (0) =1.\tag{2.1d}
\end{gather}
Here of course, and throughout (in the ODE context), dots denote
differentiations with respect to the (real) time~$t$, while (in the
ODE context) we always use primes to denote differentiations with
respect to the (complex) independent variable $\tau$ (hence
$\dot{\varphi} = \exp ({\rm i} \omega t) \varphi'$, see (2.1b)).

Consider the nonlinear first-order ODE
\begin{gather}
\varphi' = \alpha \varphi^{p/q} , \qquad \varphi \equiv \varphi
(\tau).
\end{gather}
Here $\alpha$ is an arbitrary (possibly complex) constant $(\alpha
\not= 0)$, and $p$, $q$ are two arbitrary integers. Without loss of
generality we assume $q$ to be positive, $q > 0$, and $p$, $q$ to be
coprime (namely, their decompositions into products of
primes contains no common factor). We also hereafter assume, for
simplicity, $p \not= q$, to exclude the trivial, linear, case, see~(2.5).

The initial-value problem for (2.5) is solved by the following formula:
\begin{subequations}
\begin{gather}
\varphi (\tau) = \varphi (0) [1 - (\tau / \tau_{b})]^{q/(q-p)},\\
\tau_{b} =[q/ (p-q)] \alpha^{-1} [\varphi (0)]^{(q-p)/q}.
\end{gather}
\end{subequations}
This solution, considered as a function of $\tau$, is holomorphic
inside the disk $D$ centered, in the complex $\tau$-plane, at the
origin, $\tau =0$, and having radius $\rho$,
\begin{gather}
\rho = \vert \tau_{b} \vert;
\end{gather}
indeed $\varphi (\tau)$ has (only) two (algebraic) branch points, one
at $\tau = \infty$, the other at $\tau = \tau_{b}$, see (2.6b).
Hence $\varphi$, {\it considered as a function of $t$}, is {\it
periodic} with period $t_{p}$, see (2.2), if the initial datum,
$\varphi (0)$, entails $\rho > 2/\omega$, see (2.7) and (2.6b),
consistently with the above discussion. But this condition, while
sufficient to guarantee periodicity with period~$t_{p}$, is more
stringent than necessary. A more appropriate, less stringent
but nevertheless sufficient condition to insure periodicity
(with period~$t_{p}$) of $\varphi$ as a function of $t$, is that the
branch point $\tau_{b}$, see (2.6b), fall outside the circular
contour $C$ (centered at $\tau_{c} ={\rm i}/ \omega$, and of radius
$r=1/\omega$), traveled by $\tau$ when $t$ varies over one period
$t_{p}$, namely that the distance of $\tau_{b}$ from $\tau_{c}={\rm
i}/\omega$, $|\tau_{b}-\tau_{c}|$, exceed the radius $r=1/\omega$:
\begin{gather}
|\tau_{b} -{\rm i}/\omega|>1/\omega,
\end{gather}
of course with $\tau_{b}$ defined by (2.6b).

If instead the initial datum, $\varphi (0)$, entails via (2.6b) that
the inequality (2.8) gets reversed,
\begin{gather}
\vert \tau_{b} - {\rm i} /\omega \vert < 1 / \omega,
\end{gather}
then $\varphi$, considered as a function of $t$, is also periodic, but
now with period
\begin{gather}
\tilde{t}_{p} = \vert q - p \vert t_{p};
\end{gather}
indeed in this case, (2.9), as $t$ varies over one period $t_{p}$,
$ \tau$ traverses once (before recovering its original value) the
branch cut that runs from $\tau = \tau_{b}$ to $\tau = \infty$,
hence the function $\varphi (\tau)$ picks up a phase factor $\exp
[2 \pi {\rm i} q/ (q-p)]$, see~(2.6a); hence after $\vert q - p
\vert$ such crossings the function $\varphi$ recovers its original
value.

The initial data, $\varphi (0)$, that satisfy via (2.6b) the {\it
inequalities} (2.8) respectively (2.9) are separated by a
(lower-dimensional) set of initial data, $\varphi (0)$, that satisfy
via (2.6b) the {\it equality}
\begin{gather}
\vert \tau_{b} - {\rm i} /\omega \vert = 1/\omega.
\end{gather}
This equality entails that the time $t_{b}$, defined mod ($t_{p}$)
by setting (see (2.1a) and (2.6b))
\begin{gather}
\tau_{b}= [\exp ({\rm i} \omega t_{b}) -1] /({\rm i}\omega),
\end{gather}
is {\it real}, since this formula, (2.12), via (2.11) entails
\begin{gather}
\vert \exp ({\rm i} \omega t_{b}) \vert =1.
\end{gather}
Hence for the set of initial data $\varphi (0)$ that satisfy (2.11)
(with $\tau_{b}$ defined by (2.6b)) the solution $\varphi$ diverges
(if $p > q$) or vanishes (if $p<q$) at the real time $t = t_{b} \;\,{\rm mod}\,
(t_{b})$; this signifies that, at these times, the ODE (2.5)
becomes {\it singular}.

Before proceeding further, let us note that much of the above
treatment would apply even if the {\it rational} exponent $p/q$ in
(2.5) were replaced by an arbitrary {\it real}, or even {\it
complex}, number: in the formulas above, this is most easily realized
by setting $q=1$ and letting $p$ be an arbitrary, even complex,
number. The conclusion about the {\it periodicity}, with period
$t_{p}$, of the solution $\varphi$ (considered as a function of $t$),
would still stand, whenever the intial datum $\varphi (0)$ entails,
via (2.6b), validity of the inequality (2.8); on the other hand, if
the initial datum, $\varphi (0)$, entails via (2.6b) the reversed
inequality (2.9), then $\varphi$, considered as a function of $t$,
would cease altogether to be periodic (unless $p$ is rational). The
fact that the initial data which yield (2.8) respectively (2.9) are
separated by a (lower dimensional) set of initial data, $\varphi
(0)$, satisfying (2.11) via (2.6b), data which yield solutions that
become singular at $t = t_{b}$, see (2.12), would also stand. However,
for the results that now follow, the restriction to a {\it rational}
exponent $p/q$ in the right-hand-side of (2.5) plays an essential
role: without this restriction no solution of the evolution equations
(2.17), see below, is periodic.

Let us now obtain the ODE, with the real time $t$ as independent
variable, satisfied by the (complex) function $w(t)$, see (2.4).
Clearly this definition, (2.4), entails, via (2.1b),
\begin{gather}
\dot{w} = {\rm i} \lambda \omega w + \exp [{\rm i} (1+\lambda) \omega t]
\varphi',
\end{gather}
hence, via (2.5) and (2.4),
\begin{gather}
\dot{w} - {\rm i} \lambda \omega w = \alpha \exp \{{\rm i} [1 + \lambda
(q -p)/ q] \omega t\} w^{q/p}.
\end{gather}
It is now clearly convenient to set
\begin{subequations}
\begin{gather}
\lambda = q /(p -q),\\
\Omega = \lambda \omega = [q / (p -q)] \omega,
\end{gather}
\end{subequations}
so that the evolution equation satisfied by $w \equiv w (t)$ take the
neater ({\it autonomous}!) form (1.1):
\begin{gather}
\dot{w} - {\rm i} \Omega w = \alpha w^{p/q}.
\end{gather}

The solution of the initial-value problem for this equation then
follows from (2.4), (2.6), (2.1a) and (2.16):
\begin{gather}
 w(t) = w(0) \exp ({\rm i}\Omega t) \{1\! - \!\alpha [w(0)]^{(p-q)/q} [\exp
\{{\rm i} [(p-q)/q] \Omega t\} \!-\!1]/({\rm
i}\Omega)\}^{q/(q-p)}.\!\!
\end{gather}
And it is clear, by inspection or from the previous discussion (also
keeping in mind (2.4) and (2.16)), that the set of initial data
$w(0)$ such that there hold the {\it inequality}
\begin{gather}
\left\vert \alpha^{-1} [w(0)]^{(q-p)/q} - {\rm i}/ \Omega \right\vert >
1/\Omega
\end{gather}
yield solutions (2.18) which are periodic with the period $T_{1}$,
equal to the minimum common integer multiple among $t_{p} =
2\pi/\omega = [q/(p-q)] 2 \pi/\Omega = \vert q / (p-q) \vert T$
and $T$, see (2.2) and~(1.47); while the set of initial data
$w(0)$ such that the reversed {\it inequality} hold,
\begin{gather}
\left\vert \alpha^{-1} [w(0)]^{(q-p)/q} - {\rm i}/\Omega \right\vert <
1/\Omega,
\end{gather}
yield solutions (2.18) which are periodic with period $T_{2} = q
T$, which is indeed the minimum common integer multiple among
$\tilde{t}_{p} = \vert q - p \vert 2 \pi / \omega = q 2 \pi /
\vert \Omega \vert = q T$ and $T$, see (2.10) and~(1.47); and
these two sets of intial data are separated by (the lower
dimensional set of) those data that satisfy the {\it equality}
\begin{gather}
\left\vert \alpha^{-1} [w (0)] ^{(q-p)/q} - {\rm i} / \Omega \right\vert = 1 /
\Omega,
\end{gather}
namely by the initial data $w (0)$ which yield solutions $w (t)$ that
become singular at $t = t_{b}\;\, {\rm mod}\, (t_{p})$, with $t_{b}$ given by
(2.12) and (2.6b) (with $\varphi (0)$ replaced by $w (0)$, since
these two quantities coincide, see (2.4) and (2.1c)), and with
$t_{p}$ given by (2.2) (with (2.16b)).

The ODE (2.17) is the modified version of (2.5), which indeed
features a lot of {\it periodic} solutions: in fact, as we just saw,
{\it all} its {\it nonsingular} solutions are {\it periodic} (with
periods~$T_{1}$ or $T_{2}$). Let us re-emphasize that, in (2.17), the
independent variable~$t$ (``time'') is {\it real}, while the dependent
variable, $w \equiv w (t)$, is {\it complex}, as entailed by the fact
that $\Omega$ is {\it real}; since we assumed $\omega$ and $q$ to be
{\it positive}, $\Omega$ has the same sign as $p-q$, see (2.16b); but
this is not a significant restriction on (2.17), indeed the
assumptions that $\omega$ and~$q$ be {\it positive} are quite
unessential (it is of course essential that neither $\omega$ nor $q$
vanish!); these assumptions were made merely to marginally simplify
our presentation. We can also assume the arbitrary (nonvanishing!)
constant $\alpha$ to be {\it complex}; while the exponent $p/q$ is of
course {\it rational}. Other avatars of this ODE, (2.17), are
displayed in the preceding Section~1, see (1.2), (1.3), (1.4).

Let us now proceed and show how the trick works for PDEs, beginning
from a very simple, solvable, example.

Consider the (``shock'') evolution PDE
\begin{gather}
\varphi_{\tau}= \alpha \varphi_{x} \varphi^{p/q}, \qquad \varphi \equiv
\varphi (x, \tau).
\end{gather}
Here $\alpha$ is again an arbitrary (possibly complex) constant and
$p/q$ is a (nonvanishing) {\it rational} number ($q$ a {\it positive
integer}, $q>0$; $p$ a {\it nonvanishing integer}, $p \neq 0$; $p$
and $q$ coprime). We now set (see~(2.4))
\begin{gather}
w (x, t)= \exp({\rm i} \lambda \omega t) \varphi [x, \tau(t)],
\end{gather}
of course with $\tau(t)$ defined by (2.1a) and with $\lambda$ a
parameter to be determined (see below). This entails, via (2.1b),
\begin{subequations}
\begin{gather}
w_{t}={\rm i} \lambda \omega w + \exp [ {\rm i} (1+ \lambda) \omega
t] \varphi_{\tau},
\end{gather}
hence, via (2.22) and (2.23),
\begin{gather}
w_{t}={\rm i} \lambda \omega w +\alpha  \exp [{\rm i} (1- \lambda p/q) \omega t] w_{x} w^{p/q}.
\end{gather}
\end{subequations}
This suggests setting
\begin{subequations}
\begin{gather}
\lambda = q/p,\\
\Omega = \lambda \omega = (q/p) \omega,
\end{gather}
\end{subequations}
so that the evolution PDE satisfied by $w(x,t)$ take the neat form
\begin{gather}
w_{t} - {\rm i} \Omega w = \alpha w_{x} w^{p/q}.
\end{gather}

This is the evolution PDE, see (1.35), we expect shall possess a lot of periodic
solutions. Indeed, as can be easily verified, the solution of the
initial-value (or ``Cauchy'') problem for (2.26) is given by the
definition (2.23) with (2.1a) and with $\varphi (x,\tau)$ being the
root of the (nondifferential) equation
\begin{gather}
\varphi = w_{0}\left(x+ \alpha \tau \varphi^{p/q}\right),
\end{gather}
where $w_{0}(x)$ provides the initial condition for (2.22) and (2.26)
(see (2.23) and (2.1c)):
\begin{gather}
w(x,0) = \varphi(x,0)=  w_{0}(x).
\end{gather}
But the equation (2.27) generally has several (indeed, quite
possibly, an infinity of) complex roots. Which one should be chosen?
Of course, the one that develops by continuity from (2.28) at $\tau =
t=0$, as the time evolution unfolds. Indeed it is clear that, if
$w_{0}(x)$ is an {\it analytic} function of its argument, $x$, with no
singularities for {\it real}~$x$, and if moreover the ({\it positive})
quantity
\begin{gather}
\Delta (x) = \vert \alpha | \, \rho\, |w_{0}(x)|^{p/q}
\end{gather}
is sufficiently {\it small} compared to the distance, say $d(x)$, of
the singularity of $w_{0}(z)$ in the complex $z$-plane closest to the
real point $x$,
\begin{gather}
 \Delta (x) \ll d(x),
\end{gather}
then the solution $\varphi$ of
(2.27) that flows from (2.28) can be obtained from (2.27) by
iteration, namely by setting
\begin{subequations}
\begin{gather}
\varphi_{0} (x,\tau) = w_{0} (x),\\
\varphi_{n} (x,\tau) = w_{0} \left(x + \alpha \tau [\varphi_{n-1}
(x,\tau)]^{p/q}\right), \qquad n=1,2,3,\ldots
\end{gather}
\end{subequations}
so that the sequence $\varphi_{n} (x,\tau)$ converge to
$\varphi(x,\tau)$ as $n\rightarrow \infty$,
\begin{gather}
\varphi(x,\tau) = \lim_{n\rightarrow\infty}   [\varphi_{n}
(x,\tau)]
\end{gather}
for all values of the {\it complex} variable $\tau$ such that $\vert
\tau \vert \leq \rho$, see (2.29). This entails that, if the
condition (2.30)
holds for {\it all} ({\it real}) values of $x$ (we are implicitly
assuming here that one considers the PDE (2.17) for {\it all real}
values of the space variable $x$), then $\varphi(x,\tau)$ is
holomorphic in $\tau$ in a disk $D$, centered at $\tau=0$, whose
radius $\rho$ can be made arbitrarily large by making an appropriate
choice of the initial datum, see (2.28) and (2.29). This entails that there
always is a set of initial data, $w_{0} (x)$, the modulus of which is,
for all real values of $x$, sufficiently small (if $p/q > 0$) or
sufficiently large (if $p/q < 0$) to entail that $\rho$ exceed
$2/\omega$ 
$$\rho  > 2/\omega\,.\eqno(2.33)$$
\noindent
But we know from the above discussion
that, whenever this happens, $\varphi$, considered as a function of
$t$, is periodic with period $t_{p}$, see (2.2). Hence the solution
$w(x,t)$ of (2.26) with such initial datum, see (2.28), is
periodic in $t$ with period
\setcounter{equation}{33}
\begin{gather}
\tilde{T} = \vert p \vert t_{p} = q T.
\end{gather}
see (2.23), (2.25), (2.2) and (1.47). This fulfills our expectation
that the nonlinear evolution PDE (2.26) possess (a lot of) periodic
solutions; note that the set of initial data that, according to the above
discussion, yield {\it completely periodic} solutions
$w(x,t)$ (namely, solutions periodic in $t$ with period
$\tilde{T}$, see (2.34), for {\it all real} values of $x$), have
nonvanishing measure among all initial data, $w_{0}(x)$, since the
fundamental restriction characterizing them is validity of the
{\it inequality} (2.30) for all {\it real} values of $x$, with the
quantity $\rho$ in the right-hand side of the definition (2.29) of
$\Delta(x)$ replaced by $2/\omega$.

While this discussion is sufficient to justify our claim that the nonlinear
evolution PDE (2.26)
possess a lot of completely periodic solutions, more can be said if
attention is restricted to special sets of initial data. Assume for
instance that the initial datum, $w_{0} (x)$, is a rational function:
\begin{gather}
w_{0}(x) = P(x) / Q (x),
\end{gather}
where $P(x)$, $Q(x)$ are two polynomials (and of course $Q(x)$ has no
real zeros, $Q(x) \neq 0$ for real $x$, so that $w_{0} (x)$, see
(2.35), is {\it nonsingular} for real $x$). Then (2.27) becomes
\begin{gather}
\varphi\, Q \left(x+ \alpha \tau \varphi^{p/q}\right) = P \left(x + \alpha \tau
\varphi^{p/q}\right),
\end{gather}
which is a polynomial equation for the unknown quantity
\begin{gather}
y=\varphi^{{\rm sign}(p)/q}.
\end{gather}
The degree $N$ of this polynomial equation in $y$ is clearly the
largest of the two positive integers $q+N^{(Q)} \vert p\vert$,
$N^{(P)}\vert p \vert$ if $p$ is {\it positive}, and instead the
largest of the two positive integers $N^{(Q)}\vert p \vert$, $q + N^{(P)}
\vert p \vert$ if $p$ is negative, with $N^{(P)}$ respectively
$N^{(Q)}$ the degree of the polynomial $P(x)$ respectively $Q(x)$ (of
course in order that $w_{0}(x)$ be localized in space, namely that it
vanish as $x \rightarrow \pm \infty$, the degree of $Q(x)$ should
{\it exceed} the degree of $P(x)$, $N^{(Q)} > N^{(P)}$). The
coefficients of the polynomial equation (2.36) in $y$, see (2.37), are
of course {\it periodic} in $t$ with period $t_{p}$, see (2.1a) and
(2.2). Hence all its zeros will also be {\it periodic}, with the same
period $t_{p}$, or a period which is a (finite integer) multiple of
$t_{p}$ if during the time evolution the zeros get reshuffled. Hence
we can conclude that {\it all} {\it nonsingular} solutions of the nonlinear
evolution PDE (2.26) evolving from a {\it rational} initial datum
$w_{0}(x)$ will be {\it completely periodic} (note that in this
rational case there is no restriction on the overall size of the
initial datum $w(x)$).

The third, and last, example we discuss (in this Section~2 meant to
illustrate how the trick works) takes as starting point the
(``Burgers'') evolution PDE
\begin{gather}
\varphi_{\tau} - \beta \varphi_{\xi\xi} = 2 \alpha \varphi_{\xi}
\varphi,\qquad \varphi \equiv \varphi (\xi,\tau).
\end{gather}

The arbitrary (possibly complex, but nonvanishing) constants
$\alpha$,  $\beta$ could be rescaled away, but we prefer to keep them
(as well as the standard factor 2 in the right-hand side). Here of
course subscripted variables denote partial differentiations.

Note that we have introduced a new ``space'' variable, $\xi$. The
motivation for doing so is that, in addition to the change of
variable (2.1a) (from the {\it complex} independent variable~$\tau$
to the {\it real} variable $t$, ``time''), and to the change of
dependent variable analogous to (2.23) and (2.4),
\begin{gather}
w (x,t) = \exp ({\rm i} \lambda \omega t) \varphi [\xi (t), \tau (t)],
\end{gather}
it is now convenient to introduce as well a (time-dependent)
rescaling of the ``space'' independent variable, by setting
\begin{gather}
\xi = x \exp ({\rm i} \mu \omega t).
\end{gather}
As we will immediately see this has the advantage to yield, via a
convenient choice of the two numbers $\lambda$ and $\mu$, see (2.39)
and (2.40), quite a neat form for the nonlinear evolution PDE
satisfied by $w (x,t)$ (but there is also a drawback, see below). Note
that (2.39) and (2.40) entail, via (2.1c),
\begin{gather}
\varphi (x,0) = w (x,0).
\end{gather}

Clearly (2.39) also entails, via (2.1b) and (2.40),
\begin{subequations}
\begin{gather}
w_{t} = {\rm i} \lambda\,\omega w + {\rm i} \mu \omega x \exp [{\rm i} (\mu + \lambda) \omega
t] \varphi_{\xi} + \exp [{\rm i}(1+\lambda) \omega
t]\varphi_{\tau},
\end{gather}
while (2.39) with (2.40) entail
\begin{gather}
w_{x} = \exp [{\rm i}(\mu + \lambda ) \omega t] \varphi_{\xi},\\
w_{xx}= \exp [{\rm i}(2\mu + \lambda ) \omega t] \varphi_{\xi \xi}.
\end{gather}
\end{subequations}
Hence the evolution PDE (2.38) satisfied by $\varphi (\xi , \tau)$
translates into the following evolution PDE satisfied by $w (x,t)$:
\begin{gather}
w_{t} - {\rm i} \lambda \omega w - {\rm i} \mu \omega x w_{x} - \beta
\exp [{\rm i}
(1 - 2\mu) \omega t] w_{xx} = 2 \alpha \exp [{\rm i} (1 - \mu - \lambda)
\omega t] w_{x} w.
\end{gather}
It is therefore natural to set
\begin{subequations}
\begin{gather}
\lambda = \mu = 1/2,\\
\omega = 2 \Omega,
\end{gather}
\end{subequations}
so that (2.43) take the following neat form (see (1.38)):
\begin{subequations}
\begin{gather}
w_{t} - {\rm i} \Omega (w + x w_{x}) - \beta w_{xx} = 2 \alpha w_{x} w,
\end{gather}
or equivalently
\begin{gather}
w_{t} = \left(\beta w_{x} + {\rm i} \Omega x w + \alpha w^{2}\right)_x.
\end{gather}
\end{subequations}

This is the evolution PDE we expect shall possess a lot of
solutions, $w \equiv w (x,t)$, {\it completely periodic} in the {\it
real} ``time'' $t$. This hunch is entailed by the above derivation of
this evolution PDE (from the evolution PDE (2.38), via the changes of
dependent and independent variables (2.39), (2.1a) and (2.40), with
(2.44)), and it is indeed demonstrated by the solutions we exhibit
below. Note however that this evolution PDE, (2.45), features an
explicit dependence on the space coordinate, $x$; hence it lacks
translation invariance. This explicit dependence originates from the
change of (independent, space) variable (2.40), which is on the other
hand instrumental to eliminate, via the assignment (2.44a), any
explicit time dependence from (2.45): see (2.43), and note that
forsaking the change of (space, independent) variable (2.40) from $\xi$
to $x$ amounts to setting $\mu =0$. But this change of variables,
(2.40) with (2.44), has another, possibly unpleasant, consequence: it
may cause the solutions $w \equiv w (x,t)$ of (2.45) (which are
generally obtained via (2.39) from the solutions of the solvable Burgers
equation (2.38)) to lose the property of (space) localization (namely
the property to vanish when the {\it real} variable diverges, $x
\rightarrow \pm \infty$), and it may also cause the solutions, $w
\equiv w (x,t)$, to become {\it singular} for some finite {\it real}
values of the variables $x$ and~$t$: note that, as $x$ varies over the
real axis, $- \infty < x < \infty$, and~$t$ varies over one period,
say $0 \leq t < t_{p}$, see (2.2), the variable $\xi$, see (2.40),
sweeps (in fact, twice) the {\it entire complex} $\xi$-plane; while
of course $\varphi (\xi , \tau)$ cannot be singularity-free in the
{\it entire} complex $\xi$-plane, if it does depend at all upon the
variable $\xi$; hence the singularities of $\varphi (\xi , \tau)$ in
the complex $\xi$-plane are likely to show up, via (2.39), as
singularites of $w (x,t)$ for {\it real} $x$ and $t$.
This phenomenon, which might severely reduce the applicability of
(2.45) (and of all the other evolution PDEs exhibited in this paper
for the derivation of which the change of variable (2.40) was
instrumental; but note that the evolution PDE discussed above, see
(2.26), does {\it not} belong to this class), is indeed exhibited by
many of the solutions of (2.45) displayed below; however, not by {\it all}
these solutions, as we now show.

The solutions $w \equiv w (x,t)$ of (2.45) can be obtained via (2.39)
from the solutions $\varphi \equiv \varphi (\xi , \tau)$ of the
Burgers equation (2.38). These solutions, $\varphi \equiv \varphi
(\xi, \tau)$, can themselves be obtained via the Cole--Hopf
transformation,
\begin{gather}
\varphi (\xi , \tau) = \alpha^{-1} \psi_{\xi} (\xi , \tau) / \psi
(\xi , \tau),
\end{gather}
from the solutions of the {\it linear} (``heat'', or
``Schr\"odinger'') PDE
\begin{gather}
\psi_{\tau} = \beta \psi_{\xi \xi} , \psi \equiv \psi (\xi, \tau).
\end{gather}
In this manner one can easily manufacture large classes of
solutions $\varphi (\xi , \tau)$ of the Burgers equation (2.38),
hence, via (2.39), large classes of solutions $w (x,t)$ of the
evolution PDE (2.45); and one can also generally solve, in this
manner, the initial-value problem for (2.45). For the reasons explained
above, one then expects many of these solutions $w (x,t)$ to depend
{\it periodically} on the real time $t$, but also to become singular
for some {\it real} values of the ``space'' and ``time'' variables
$x$ and $t$. But there also exist solutions of (2.45) which are
{\it nonsingular} for {\it all real} values of $x$ and $t$, and which
are moreover {\it localized} in $x$ (namely, such that $w ( \pm
\infty , t) =0$) and of course {\it periodic} in~$t$.

These statements are validated by the following two explicit examples.

The solution $w (x,t)$ of (2.45) that corresponds to the solution
\begin{gather}
\psi (\xi , \tau) = A - \exp \left(k \xi + \beta k^{2} \tau\right)
\end{gather}
of the {\it linear} PDE (2.47), reads
\begin{subequations}
\begin{gather}
w (x,t) =- (k B /\alpha) \exp ({\rm i} \Omega t)\nonumber\\
\qquad {}\times  \left\{B - \exp \left[- kx \exp ({\rm i} \Omega t) + \left[{\rm i} \beta k^{2} /
(2 \Omega)\right] \exp (2 {\rm i} \Omega t)\right]\right\}^{-1},\\
B = A^{-1} \exp \left[-{\rm i} \beta k^{2} / (2 \Omega)\right].
\end{gather}
\end{subequations}
Here $k$ and $B$ (or, equivalently, $k$ and $A$, see (2.49b)) are two {\it
arbitrary} (possibly complex) constants. This solution is obviously
{\it periodic} in $t$ with period $T$, see (1.47); but (as it can be
easily shown, see the Appendix) it {\it always} becomes singular at some
{\it real} value of $x$ and~$t$.

Another interesting class of solutions of (2.45) are characterized by
the property to be, for all time, {\it rational} in $x$. They are
easily obtained, via (2.39), (2.40) and (2.1a), from (2.46) with
$\psi (\xi, \tau)$ a solution of the linear PDE (2.47) which is, for
all time, polynomial in $\xi$ (say, of degree $N$):
\begin{gather}
\psi (\xi , \tau) = \sum^{N}_{m=0} \gamma_{m} (\tau) \xi^{N -m}.
\end{gather}
It is easily seen that such polynomial solutions of (2.47) exist and
that their coefficients are characterized by the ODEs
\begin{gather}
{\gamma}'_{0} = {\gamma}'_{1} = 0, \qquad {\gamma}'_{m} =\beta
(N-m +2) (N-m+1)\gamma_{m-2} , \qquad m=2,\ldots,N,
\end{gather}
which could be easily solved in explicit form.

These rational solutions $w (x,t)$
 of (2.45) are clearly {\it periodic} in $t$
with period $T$ (see (2.39), (2.40), (2.1a), (2.44) and (1.47)), and they
are, for all time, {\it spatially localized}, namely they vanish as
$x \rightarrow \pm \infty $, $w ( \pm \infty, t)=0$ (albeit the
vanishing is slow, being generally proportional to $\vert x
\vert^{-1}$ as $x \rightarrow \pm \infty$).
But is there any such solution which remains {\it nonsingular} for
all {\it real} values of $x$ and $t$?

To gain some insight about this question, we note that {\it all}
these rational solutions $w (x,t)$ also admit the following
representation:
\begin{gather}
w(x,t) = (\beta /\alpha) \sum^{N}_{n=1} [x - z_{n} (t)]^{-1},
\end{gather}
with the $N$ {\it complex} quantities $z_{n} (t)$ characterized by
arbitrary ``initial conditions'', $z_{n} (0) = z_{n}^{(0)}$, and by the
following system of first-order ODEs:
\begin{gather}
\dot{z}_{n} + {\rm i} \Omega z_{n} = - 2\beta \sum^{N}_{m=1 , m \not= n}
(z_{n} - z_{m})^{-1}.
\end{gather}
As it is well known [3], and as it is indeed implied by the above
treatment, these $N$ quantities $z_{n} (t)$ are simply related to
the $N$ zeros of the time-dependent polynomial (2.50); they are of
course all periodic functions of time, with period at most
$\tilde{T} = T \cdot N!$.

The issue about the existence of {\it nonsingular} solutions $w (x,t)$
(rational in $x$) of (2.45) hinges now on the question whether, for
given $\beta$ and ({\it real}!) $\Omega$, and for some given $N$,
there exist initial (of course {\it complex}) values $z_{n} (0)$ such
that {\it none} of the quantities $z_{n} (t)$ is {\it real} for any
real time. The structure of (2.52) might seem to suggest a negative
reply; note for instance that the center-of-mass of the quantities $z_{n}
(t)$,
\begin{gather}
Z(t) = N^{-1} \sum^{N}_{n=1} z_{n} (t),
\end{gather}
evolves in time according to the equation
$$
\dot{Z} + {\rm i} \Omega Z =0,
\eqno(2.55a)$$
hence according to the law
\setcounter{equation}{54}
$$
Z (t) = Z (0) \exp (-{\rm i} \Omega t),
\eqno(2.55b)$$
which clearly entails that $Z (t)$, whatever the choice of the initial
datum $Z(0)$, at some ({\it real}) time does indeed become {\it real}.

But the pessimism about the existence of {\it nonsingular} solutions $w
(x,t)$ (rational in $x$) of (2.45) turns out to be unwarranted: there
do exist such {\it nonsingular} solutions, as shown by the following
example (for $N=2$). Indeed such a solution $w (x,t)$ of (2.45) reads,
as can be easily verified,
\begin{gather}
w (x,t) = 2 [x - a \exp (- {\rm i} \Omega t)]
\left\{ \left[x -a \exp (-{\rm i} \Omega t)
\right]^{2}\! - {\rm i} (\beta /\Omega) -b \exp (-2 {\rm i} \Omega t)\right\}^{-1}\!\!,\!\!
\end{gather}
with $a$ and $b$ two arbitrary (complex) constants; and we show in
the Appendix that a condition {\it sufficient} to guarantee that this
solutions, (2.56), be {\it nonsingular} for {\it all real} values of
$x$ and $t$, is validity of the {\it inequality}
\begin{gather}
\vert {\rm Re}\, (\beta) / \Omega \vert > \vert b \vert + 4 \vert a
\vert^{2}.
\end{gather}

Note that, if $a=0$, this solution $w (x,t)$, see (2.56), of the
evolution PDE (2.45) is an {\it odd} function of $x$, and it is {\it
completely periodic} in t with period $T/2$, see (1.47); if $a \not= 0$,
this solution $w (x,t)$, see (2.56), has no definite parity as a
function of $x$, and it is {\it completely periodic} in $t$ with
period $T$, see (1.47).

\section{Derivation of evolution equations featuring\\
periodic solutions}

In the preceding Section 2 we have discussed in some detail, using
three representative examples, the trick whereby evolution equations
can be modified so that the modified versions (have
perhaps a neat form and) generally feature a lot of, and in some
cases only, {\it periodic} solutions. In this Section~3 we apply this
trick to other evolution equations, and we thereby obtain all the
evolution equations exhibited in the introductory Section~1. Our
treatment below is however somewhat terser than in the preceding
Section~2.

We treat firstly ODEs,  then PDEs in $1+1 =2$ variables (space and
time), and finally one example in $1+2=3$ variables.

\medskip

{\bfseries\itshape Derivation of (1.5), (1.6) and (1.7.).} {\it
All} solutions of the first Painlev\'e ODE
\begin{gather}
\varphi'' = \alpha \varphi^{2} + \beta \tau, \qquad \varphi \equiv
\varphi (\tau),\end{gather}
 are {\it meromorphic} functions of the
complex independent variable $\tau$ (see, for instance, [4]). The
two arbitrary (complex) constants $\alpha$, $\beta$ could be
eliminated by rescaling the dependent and independent variables,
but we prefer to keep them. Of course here and below primes denote
differentiations with respect to~$\tau$.

We now introduce the new independent variable $t$ (``time'') by
setting \begin{gather} \tau \equiv \tau (t) = - ({\rm i} / \omega)
\exp ({\rm i} \omega t).\end{gather}
 Here and throughout $\omega$
is a {\it positive} number, $\omega > 0$, to which we associate
the period $t_{p}$, see (2.2). Note the analogy, but also the
difference, among (3.2) and (2.1a).

We moreover introduce the new (complex) dependent variable $w
\equiv w (t)$ via a relation closely analogous to (2.4):
\begin{gather}
w (t) = \exp ({\rm i} \lambda \omega t) \,\varphi [ \tau (t)],
\end{gather} where $\lambda$ is a number to be chosen
conveniently, see below. Clearly these formulas entail that, if
$\lambda$ is a {\it rational} number, the property of $\varphi
(\tau)$ to be {\it meromorphic} entails that $w (t)$ is {\it
periodic} in $t$ with a period which is an {\it integer} multiple
of $t_{p}$, see (2.2).

On the other hand, via (3.2) and (3.3), the evolution ODE (3.1) yields
\begin{subequations}
\begin{gather}
\ddot{w} = {\rm i} (2 \lambda +1) \omega \dot{w} + \lambda
(\lambda +1) \omega^{2} w + \alpha \exp [{\rm i} (2 - \lambda)
\omega t] w^{2} + \gamma \exp [{\rm i} (\lambda +3) \omega
t],\end{gather} where we set \begin{gather} \beta = {\rm i} \omega
\gamma.
\end{gather}
\end{subequations}

Of course here and below dots denote differentiations with respect
to the ({\it real}) time~$t$.

Three choices of $\lambda$ appear of special interest: $\lambda =
- 1/2$, $\lambda =2$, respectively $\lambda = - 3$. The
corresponding evolution ODEs have been reported in Section~1,
see~(1.5) (where we set $\omega = 2 \Omega$), (1.6) (where we set
$\omega = \Omega$) respectively (1.7) (where we also set $\omega =
\Omega$). The (extremely simple!) periodicity property of their
solutions, as implied by (3.2) and~(3.3), are mentioned there.
Note the difference among the transformations (2.1a) and~(3.2),
including the fact that, for $\omega \rightarrow 0$, (2.1a)
becomes the identity $\tau = t$ while (3.2) has no limit. This
accounts for the fact that, for $\Omega = 0$, (1.5) does $not$
reproduce the Painlev\'e equation (3.1). Analogous remarks apply
to the analogous cases treated below.

\medskip

{\bfseries\itshape Derivation of (1.8).} The treatment is closely
analogous to that described immediately above, except that one now
starts from the second Painlev\'e ODE~[4], \begin{gather}\varphi''
= \alpha \varphi^{3} + \beta \tau \varphi + \delta,\qquad \varphi
\equiv \varphi (\tau),\end{gather} rather than from the first, see
(3.1). Here again we keep some constants that could be removed by
appropriate rescalings. One then uses (3.2) and (3.3), and obtains
thereby \begin{subequations} \begin{gather} \ddot{w} = {\rm i} (2
\lambda +1) \omega \dot{w} + \lambda (\lambda +1) \omega^{2} w
\nonumber\\ \qquad {}+ \alpha \exp [2 {\rm i} (1 - \lambda) \omega
t]w^{3} + \gamma \exp (3 {\rm i} \omega t) w + \delta \exp [{\rm
i} (2 + \lambda) \omega t]
\end{gather}
with \begin{gather} \gamma = -{\rm i} \beta / \omega. \end{gather}
\end{subequations}
This suggests setting $\lambda =1$, $\omega = \Omega$, and one
thereby gets (1.8).

Of course the fact [4] that all solutions of the second Painlev\'e
ODE, (3.5), are meromorphic functions of $\tau$ entails, via (3.2)
and (3.3) (with $\lambda = 1$, $\omega = \Omega$) that {\it all
nonsingular} solutions of (1.8) are periodic functions of the
({\it real}) time $t$, with period $T$, see (1.47).

\medskip

{\bfseries\itshape Derivation of (1.13).} One starts from
\begin{gather}
\varphi'' = \alpha \varphi^{2} \end{gather} and sets
\begin{subequations}
\begin{gather}
w(t) =\exp ({\rm i} \Omega t) \varphi (\tau) \end{gather}
with
$\tau$ given by (2.1a) and
\begin{gather}
\Omega = 2 \omega \end{gather} \end{subequations} (this
corresponds to (3.3) with $\lambda =2$). This entails that the
{\it general} solution of (1.13) reads
\begin{gather}
w(t) = (6 / \alpha) \exp ({\rm i} \Omega t) {\cal{P}}[-(2 {\rm
i}/\Omega) \exp ({\rm i} \Omega t/2)+ \beta; \; 0, g_{3}],
\end{gather}
 with $\beta$ and $g_{3}$ arbitrary (complex)
constants, and \begin{gather} {\cal{P}} (z;g_{2},g_{3}) \equiv
{\cal{P}} (z \vert \omega, \omega') \end{gather}
 the Weierstrass
doubly-periodic elliptic function which satisfies the ODE
\begin{gather}
{\cal{P''}} (z) = 6 {\cal{P}}^{2} (z) - g_{2}/2.
\end{gather} Note the (standard) notation (3.10), which
entails that $\omega$ and $\omega'$ are the two semiperiods of the
Weierstrass function (note that the quantity $\omega$ in the
right-hand side of (3.10) has nothing to do with the ($real$)
quantity $\omega$ in (2.1a) and (3.8b)). Also note that a
necessary and sufficient condition to guarantee that the solution
(3.9) be $nonsingular$ for all ($real$) values of the time $t$ is
the inequality \begin{gather} \vert \beta \vert \neq \vert 2 /
\Omega \vert, \end{gather} and that $all$ nonsingular solutions
(3.9) are $periodic$ in $t$ with period $2T = 4 \pi / \Omega$, see
(1.47) (only in exceptional cases, namely for special values of
$\beta$ and $g_{3}$, they might be $periodic$ with period $T=2 \pi
/ \Omega$, due to the $periodicity$ of the Weierstrass function).

\medskip

{\bfseries\itshape Derivation of (1.14a), (1.14b).} One starts
from the equation \begin{gather} \varphi'' = \alpha \exp(\varphi),
\end{gather} and sets
\begin{gather}
w (t) = 2 {\rm i} \omega t + \varphi [\tau (t)] \end{gather}
 with
$\tau (t)$ given by (2.1a). This yields (1.14a) (with $\omega =
\Omega$). The position (3.14) suggests that the solutions $w
\equiv w (t)$ of (3.13) are not periodic in $t$, but that their
time-derivatives, $\dot{w} (t)$, are indeed periodic, with period
$T$, see (1.47). This expectation is confirmed by the following
expression of the {\it general} solution of (1.14a),
\begin{gather}
w(t) = a \exp ({\rm i} \Omega t) + b +2 {\rm i} \Omega
t\nonumber\\ \qquad{} - 2\log \left\{(\alpha / \Omega^{2}) \left(2
a^{2}\right)^{-1} \exp (b) + \exp [a \exp ({\rm i} \Omega
t)]\right\},\end{gather} which features the two {\it arbitrary}
constants $a$, $b$, and which entails
\begin{gather}
\dot{w} (t) = {\rm i} \Omega a \exp ({\rm i} \Omega t) + 2{\rm i} \Omega \Big[1 -
a \exp ({\rm i} \Omega t) \exp [a \exp ({\rm i} \Omega
t)]\nonumber\\ \qquad{}\times \left\{\left(\alpha / \Omega^{2}
\right) \left(2 a^{2}\right)^{-1}\exp (b) + \exp [a \exp ({\rm i}
\Omega t)]\right\}^{-1}\Big].\end{gather}
 Clearly all nonsingular
functions $\dot{w} (t)$, see (3.16), are periodic in $t$ with
period $T$, see (1.47), and it is also easily seen that a
condition sufficient to guarantee that the solution $w (t)$, see
(3.15), as well of course as its time-derivative, see (3.16), be
{\it nonsingular} for {\it all} ({\it real}) times, is validity of
the inequality \begin{gather} \vert a \vert \not= \left\vert b -
\log \left(-2 a^{2} \Omega^{2} / \alpha\right)\right \vert.
\end{gather}

To obtain (1.14b), time-differentiate (1.14a), eliminate $\alpha
\exp(w)$ using (1.14a), and then replace formally (as a notational
change) $\dot{w} (t)$ with $w(t)$. This derivation entails of
course that the {\it general} solution, $w = w(t)$, of (1.14b) is
provided by the right-hand side of (3.16); it is clearly {\it
nonsingular} if the inequality (3.17) holds, and periodic in $t$
with period $T$, see~(1.47).

\medskip

{\bfseries\itshape Derivation of (1.20), (1.21), (1.22) and (1.23).}
 Here we take as starting point the ODE
\begin{gather}
\varphi'' = \alpha (\varphi')^{p_{2}/q_{2}} \varphi^{p_{1}/q_{1}},
\qquad \varphi \equiv \varphi (\tau),
\end{gather}
where, as usual, we keep the {\it arbitrary} (nonvanishing)
constant $\alpha$ in spite of the fact that it could be easily
rescaled away, and we introduce the 4 {\it integers} $p_{1}$, $q_{1}$,
$p_{2}$, $q_{2}$ to write the two {\it rational} exponents in the
right-hand side of this ODE, (3.18). We do not dwell on the
restrictions to be imposed on these integers, which we trust are
self-evident from (1.20) as well as from the equations written below.

We now use the change of independent variable (2.1a), as well as the
change of dependent variable (3.3), the latter with
\begin{subequations}
\begin{gather}
\lambda = - q_{1}(2q_{2}-p_{2})/ (q_{1}
q_{2} - p_{1} q_{2} - p_{2} q_{1}),\\
\lambda \omega =  \Omega,
\end{gather}
\end{subequations}
and we thereby get for $w \equiv w (t)$ the ({\it autonomous}!)
ODE (1.20). Hence the solutions of this ODE, (1.20) , can be obtained,
via (2.1) and (3.3) with (3.19), from the general solution of (3.18),
which is yielded by the quadrature formula
\begin{gather}
\int^{\varphi} dx \left[x^{(p_{1}+q_{1})/q_{1}} +
a^{2}\right]^{q_{2}/(p_{2}-2q_{2})}\nonumber\\ \qquad{}= b + \{-
\alpha q_{1}(p_{2}-2q_{2})/[(p_{1}+q_{1})q_{2}]
\}^{-q_{2}/(p_{2}-2q_{2})} \tau,
\end{gather}
where $a^{2}$ and $b$ are two {\it arbitrary} constants.

We forsake here a discussion of (3.20) for an arbitrary choice of the
two rational numbers $p_{1} /q_{1}$ and $p_{2}/q_{2}$, and we limit
our consideration to three examples.

The first obtains by setting $p_{1}/ q_{1} = - 3$, $p_{2} / q_{2} =3$,
entailing $\lambda = 1$ (see 3.19a));
thereby (1.20) becomes (1.21). In this case (3.20) reads
\begin{subequations}
\begin{gather}
-\varphi^{-1} + a^{2} \varphi = b + (2/ \alpha ) \tau,
\end{gather}
entailing, say,
\begin{gather}
\varphi = \left[b +(2/ \alpha ) \tau - \left\{[b + (2/ \alpha)
\tau]^{2} + 4a^{2} \right\}^{1/2}\right]/ \left(2a^{2}\right),
\end{gather}
\end{subequations}
so that (see (2.1) and (3.3), (3.19))
\begin{subequations}
\begin{gather}
\varphi (0) = w (0) = \left[b-\left(b^{2} +
4a^{2}\right)^{1/2}\right] / \left(2a^{2}\right),
\\
\varphi' (0) = \dot{w} (0) -  {\rm i} \Omega w (0) =
\left[1-b\left(b^{2} + 4a^{2}\right)^{-1/2}\right]/
\left(\alpha a^{2}\right);
\end{gather}
\end{subequations}
and it is easily seen that the initial data, $w (0)$ and $\dot{w} (0)$,
of (1.21) are split into two sets, all of which however yield
nonsingular solutions periodic
with period $T$, see (1.47). These two sets are separated by
a (topologically nontrivial, lower dimensional) set of initial data,
characterized, see (3.21) (with (2.1a), (3.19b) and $\lambda =1$), by
the {\it equalities}
\begin{gather}
\vert 1 +  \Omega \alpha ( \pm a - {\rm i} b/2) \vert =1;
\end{gather}
and clearly these special initial data yield solutions which become
singular (namely, which are such that $\dot{w}$ diverges) at the {\it
real} times $t_{b}$, characterized by the relation
\begin{gather}
\exp ({\rm i} \Omega t_{b}) = 1 +  \Omega \alpha (\pm a - {\rm i}
b/2)
\end{gather}
(the fact that the values of $t_{b}$, as defined mod$\,(T)$ by this
equation, (3.24), are {\it real} is of course implied by (3.23)).

The second example obtains by setting
\begin{subequations}
\begin{gather}
p_{1} /q_{1} = m ,\qquad p_{2}/ q_{2} =(2n+1)/n,
\end{gather}
entailing (see (3.19a))
\begin{gather}
\lambda = -  (nm + n+1)^{-1},
\end{gather}
\end{subequations}
with $m$ a {\it nonnegative} integer and $n$ a {\it positive} integer;
thereby (1.20) becomes (1.22). In this case (3.20) reads
\begin{gather}
P_{nm+n+1} (\varphi) = b + \left\{ - [n (m+1)]^{-1} \alpha\right\}^{-n}
\tau,
\end{gather}
where $P_{nm + n+1} (\varphi)$ is a polynomial of degree $nm+n+1$ in
$\varphi$, whose coefficients are time-independent (for given $n$ and
$m$, they only depend on the constant $a$). Since $\tau$, see (2.1), is
a periodic function of $t$ with period $t_{p}$, see (2.2), clearly
the set of the $nm+n+1$ roots of (3.26) is also periodic with the
same period; but each root, if followed continuously as function of
$t$, need only be periodic with period $t_{p} \cdot (nm + n+1)!$.
Hence we may conclude that {\it all nonsingular} solutions $w(t)$ of
(1.22) are {\it periodic} with a  period $\tilde {T}$ which is at least
the minimum common integer multiple among $T$, see (1.47), and
\begin{subequations}
\begin{gather}
t_{p}= T /(nm-n+1),
\end{gather}
and at most the minimum common integer multiple among $T$, see (1.47),
and
\begin{gather}
\tilde{t}_{p} = t_{p} (nm + n+1)! = T(nm+n)!
\end{gather}
\end{subequations}
(see (3.19), (3.25) and (1.47)). The singular solutions correspond to
lower dimensional sets of initial data $w(0)$, $\dot{w}(0)$, such that,
for some {\it real} time, (at least) 2 roots of (3.26) coincide, at
which time $\dot{w} (t)$ diverges.

The third example obtains by setting
\begin{subequations}
\begin{gather}
p_{1}/q_{1} = - (2m+1), \qquad
p_{2}/q_{2} = (2n+1)/n,
\end{gather}
entailing (see (3.19a))
\begin{gather}
\lambda = (2nm-1)^{-1},
\end{gather}
\end{subequations}
with $m$ and $n$ {\it positive} integers; thereby (1.20) becomes
(1.23). In this case (3.20) reads
\begin{gather}
P_{2mn-1} (\varphi) =\left[ b+ (2mn)^{n} \alpha^{-n} \tau\right]
\varphi^{2mn-1}.
\end{gather}
Hence, by an analysis closely analogous to that given above, in this
case we also conclude that {\it all nonsingular} solutions $w(t)$
of (1.23) are {\it periodic} with a period $\tilde{T}$ which is at
least the minimum common integer multiple among $T$, see (1.47), and
\begin{subequations}
\begin{gather}
t_{p} = T /(2 nm-1)
\end{gather}
and at most the minimum common integer multiple among $T$, see (1.47),
and
\begin{gather}
\tilde{t}_{p} = t_{p} (2nm-1)! =T (2nm-2)!
\end{gather}
\end{subequations}
(see (3.19), (3.28) and (1.47)).

\medskip

{\bfseries\itshape Derivation of (1.24), (1.25) and (1.26).}
We take as starting
point the evolution ODE
\begin{gather}
\varphi'' = \alpha \varphi' \varphi + \beta \varphi^{3},\qquad \varphi
\equiv \varphi (\tau),
\end{gather}
and make the (by now standard) changes of dependent and independent
variables (3.3) and (2.1). This yields
\begin{gather}
\ddot{w}={\rm i} (2 \lambda +1) \omega \dot{w} + \lambda (\lambda
+1) \omega^{2} w \nonumber\\
\qquad{}+\alpha \exp [{\rm i} (1 - \lambda) \omega t] w
(\dot{w} - {\rm i} \lambda \omega w) + \beta \exp [2 {\rm i} (1 -
\lambda) \omega t] w^{3},
\end{gather}
and this yields (1.24) by setting
\begin{subequations}
\begin{gather}
\lambda = 1,\\
\omega = \Omega.
\end{gather}
\end{subequations}

The other two ODEs, (1.25) respectively (1.26), are clearly the two
special cases of (1.24) corresponding to $\beta = - \alpha^{2}/9$
respectively $\beta =0$. The motivation for singling out these two
cases is because, by setting
\begin{subequations}
\begin{gather}
\varphi (\tau)= \gamma \psi'(\tau) / \psi(\tau), \\
\gamma = \left[ \alpha - \left(\alpha^{2} + 8 \beta\right)^{1/2}
\right] / (2 \beta),
\end{gather}
\end{subequations}
one transforms (3.31) into
\begin{subequations}
\begin{gather}
\psi''' \psi = \eta \psi' \psi'',\\
\eta = \alpha^{2} \left\{1+6\left(\beta / \alpha^{2}\right)
- \left[ 1+8 \left( \beta /\alpha^{2}\right)\right]^{1/2}\right\} /
(2 \beta),
\end{gather}
\end{subequations}
and for $\beta =-\alpha^{2}/9$ respectively $\beta=0$ one gets from
(3.35b) $\eta =0$ respectively $\eta=1$, two values for which (3.35a)
is particularly easy to integrate.

Indeed in the first case, $\beta =-\alpha^{2}/9$, $\eta=0$ (which
entails $ \gamma =-3/ \alpha$, see (3.34b)) one gets from (3.35a)
\begin{gather}
\psi (\tau) = a + b \tau + \tau^{2}
\end{gather}
hence (see (3.34a))
\begin{gather}
\varphi (\tau) =- (3/ \alpha) ( b + 2 \tau)/ \left(a + b \tau +
\tau^{2}\right)
\end{gather}
with $a$ and $b$ two {\it arbitrary} constants. Via (3.3) (with
(3.33)) and (2.1) this entails
\begin{subequations}
\begin{gather}
w (t)= 3 {\rm i} ( \Omega / \alpha) [(2 - B) \exp({\rm i}
\Omega t) - 2 \exp (2 {\rm i} \Omega t)]\nonumber\\
\qquad{}\times [1-A-B+(B-2) \exp ({\rm i} \Omega t) + \exp (2 {\rm i} \Omega
t)]^{-1}
\end{gather}
with $A$ and $B$,
\begin{gather}
A=a \omega^{2}, \qquad B={\rm i} \omega b,
\end{gather}
\end{subequations}
two arbitrary constants. This is the general solution of (1.25);
it is clearly {\it periodic} in~$t$ with period $T$, see (1.47),
and it is {\it nonsingular} for {\it all (real)} values of $t$
provided the following {\it inequality} holds:
\begin{gather}
\left\vert 1- B/2 \pm \left[(B/2)^{2}+ A\right] ^{1/2} \right\vert
\neq 1.
\end{gather}
In the second case, $\beta =0$, the ODE (3.28) is easily integrated,
to yield
\begin{gather}
\varphi ( \tau) = (a/ \alpha) [b+\exp (a \tau)]/[b-\exp (a \tau)],
\end{gather}
with $a$ and $b$ arbitrary constants. Via (3.3) and (2.1) this entails
\begin{subequations}
\begin{gather}
w(t) = {\rm i} \Omega ( A /\alpha) \exp ({\rm i} \Omega t)
\{B+ \exp[A \exp({\rm i} \Omega
t)]\}/ \{B - \exp [  A \exp ( {\rm i} \Omega t)] \}
\end{gather}
with $A$ and $B$,
\begin{gather}
 A = - {\rm i} a / \Omega, \qquad B =b \exp (A),
\end{gather}
\end{subequations}
two {\it arbitrary} constants. This is the general solution of
(1.26); it is clearly {\it periodic} in~$t$ with period $T$, see
(1.47), and it is {\it nonsingular} provided the following {\it
inequality} holds:
\begin{gather}
\vert A \vert \neq \vert \log (B) \vert
\end{gather}
(see (3.41)).

{\bfseries\itshape Derivation of (1.29), (1.30), (1.31) and (1.32).}
Time-differentiation of (1.6) yields
\begin{gather}
\dddot{w} - 5 {\rm i} \Omega \ddot{w} - 6 \Omega^{2} \dot{w} = 2 \alpha w
\dot{w} + 5 {\rm i} \Omega \gamma \exp (5 {\rm i} \Omega t),
\end{gather}
and using again (1.6) to eliminate the last term in the right-hand
side of this equation, (3.43), one gets (1.29). Note that the
constant $\gamma$, see (1.6), does not appear in (1.29). Hence
{\it all} solutions of (1.29) also satisfy (1.6) (for some
appropriate value of $\gamma$); or, equi\-va\-lently, the solution
of the initial-value problem for (1.29) (namely, of the problem to
evaluate the solution $w(t)$ of (1.29) which corresponds to given
initial data $w(0)$, $\dot{w} (0)$, $\ddot{w} (0)$) is provided by
the solution of the initial-value problem for (1.6) (with the same
data $w(0)$, $\dot{w} (0)$, and with $\gamma = \ddot{w} (0) - 5
{\rm i} \Omega \dot{w} (0) - 6 \Omega^{2} w(0)$). Hence the
solutions of (1.29) have the same periodicity properties as the
solutions of (1.6).

As for (1.30), it is merely another avatar of (1.29). Indeed if, in
this ODE, (1.29), we set $w(t) = z(t) + c$, $\dot{c} =0$, we get
\begin{gather}
\dddot{z} - 10 {\rm i} \Omega \ddot{z} - \left(31 \Omega^{2} + 2
\alpha c\right) \dot{z} \nonumber\\ \qquad{}+ 10 {\rm i} \left(3
\Omega^{2} + \alpha c\right) \Omega z + 5 {\rm i}
\left(6\Omega^{2} + \alpha c\right) \Omega c = \alpha (2 \dot{z} -
5 {\rm i} \Omega z) z,
\end{gather}
which is a more general avatar of (1.29) than (1.30), to which it
reduces for $\alpha c = - 6 \Omega^{2}$.

The derivation of (1.31) from (1.7) is completely analogous to the
derivation given above of (1.29) from (1.6); and the comments given
above on the possibility to obtain {\it all} the solutions of (1.29)
from those of (1.6) are as well applicable now, except that the role
previously played by the constant $\gamma$ (which appears in (1.6)
but not in (1.29)) is now played by the constant $\alpha$ (which
appears in (1.7) but not in (1.31)).

The starting point to obtain (1.32) is (1.5), which we now write as
follows (by setting $\gamma = \alpha \eta$):
\begin{gather}
\ddot{w} + \Omega^{2} w = \alpha \left(w^{2} + \eta\right)
\exp (5 {\rm i} \Omega t).
\end{gather}

We now time-differentiate this ODE, (3.45), and then eliminate the
explicitly time-dependent term using again (3.45). This yields (1.32).
Again we note that the constant~$\alpha$, which appears in (3.45), has
dropped out of (1.32); hence we may again conclude that {\it all}
solutions of (1.32) also satisfy (3.45) (with an appropriate value
of $\alpha$, including $\alpha =0$), and that the solution of the
initial-value problem for (1.32) can be obtained from the solution of
the initial-value problem for (3.45) (or, equivalently, (1.5)).

\medskip

{\bfseries\itshape Derivation of (1.40), (1.41), (1.42), (1.43) and (1.44).}
Here we take as starting point the evolution PDE
\begin{gather}
\varphi_{\tau} = \beta \varphi_{\xi \xi \xi} + \alpha
\varphi_{\xi} \varphi^{p/q},\qquad
\varphi \equiv \varphi (\xi ,
\tau),
\end{gather}
with $\alpha$ and $\beta$ two {\it arbitrary} (possibly complex)
constants and $p$, $q$ two arbitrary {\it integers} $(q \not=0)$. We then
use the transformation formulas (2.1), (2.39) and (2.40). We thereby
get, for the new dependent variable $w \equiv w (x,t)$,
\begin{gather}
w_{t} = {\rm i} \lambda \omega w + {\rm i} \mu \omega x
w_{x}\nonumber\\ \qquad{} + \beta \exp [{\rm i} (1-3 \mu) \omega
t]w_{xxx} + \alpha \exp \{ {\rm i} [1 - \mu - \lambda (p/q)]
\omega t\} w_{x} w^{p/q}.
\end{gather}
Hence we set
\begin{gather}
\mu = 1/3 ,\qquad \lambda =2 p / (3q),\qquad
 \omega = 3 \Omega,
\end{gather}
and we thereby get (1.40).

The evolution PDE (1.41) is the special case of (1.40) with $q=p =1$;
for $\Omega = 0$ it reduces to the well-known ({\it integrable}!)
``Korteweg-de  Vries'' (KdV) equation.

Likewise, the evolution PDE (1.43) is the special case of (1.40) with
$q=1$, $p=2$; for $\Omega =0$ it reduces to the well-known ({\it
integrable}!) ``modified Korteweg-de Vries'' equation.

As for the two systems of two  {\it real} PDEs (1.42) respectively
(1.44), they are clearly the {\it real} avatars of (1.41)
respectively (1.43), obtained by setting $w = u +{\rm i}v$, $\alpha =
a_{1} + {\rm i} a_{2}$, $\beta = b_{1} + {\rm i} b_{2}$.

Finally, to obtain (1.45), as well as (1.46) which is just the $real$
version of the $complex$ equation (1.45) (via $w=u+ {\rm i}v$,
$\alpha = a_{1} + {\rm i} a_{2}$, $\beta= b_{1} + {\rm i} b_{2}$,
$\gamma = c_{1} + {\rm i} c_{2})$, we start from the
Kadomtsev--Petviashvili (KP) equation
\begin{gather}
(\varphi_{\tau} + \beta \varphi_{\xi\xi\xi}+\alpha
\varphi_{\xi}\varphi )_{\xi}+ \gamma \varphi_{\eta \eta}=0
\end{gather}
and we then set
\begin{subequations}
\begin{gather}
w (x,y,t) = \exp ({\rm i}\Omega t)  \varphi (\xi, \eta, \tau)
\end{gather}
with
\begin{gather}
\xi = x \exp ({\rm i} \Omega t /2),\\
\eta = y \exp ({\rm i} \Omega t),
\end{gather}
and $\tau \equiv \tau (t)$ given by (2.1a) with
\begin{gather}
\Omega = 2 \omega, \qquad \omega = \Omega/2.
\end{gather}
\end{subequations}

\appendix
\renewcommand{\theequation}{A.\arabic{equation}}
\setcounter{equation}{0}
\section*{Appendix}

In this Appendix we prove some results whose detailed treatment in the
context of the paper would have been too distracting.

The first result we prove is that the solution (2.49) of the
evolution PDE (2.45) is singular for some {\it real} value of $x$
and $t$, namely that, for {\it any} choice of the {\it complex} numbers
$\alpha$, $\beta$, $k$ and $B$, and of the {\it real} number $\Omega$,
there exist some {\it real} value of $x$ and $t$ (the latter, of
course, defined ${\rm mod}\, (T)$, see (1.47)) such that
\begin{subequations}
\begin{gather}
B = \exp \left\{- kx \exp ({\rm i} \Omega t) +
\left[{\rm i} \beta k^{2} / (2\Omega)\right]
\exp (2{\rm i} \Omega t)\right\}.
\end{gather}
Indeed this equation entails
\begin{gather}
\log (B) = - k x \exp ({\rm i} \Omega t) +
\left[{\rm i} \beta k^{2} / (2 \Omega)\right]
\exp (2{\rm i} \Omega t)
\end{gather}
namely
\begin{gather}
x = a \exp (-{\rm i} \Omega t) + b \exp ({\rm i} \Omega t),
\end{gather}
with
\begin{gather}
a=-k^{-1} \log (B), \qquad
b={\rm i} \beta k/ (2 \Omega).
\end{gather}
We  now show that, for any arbitrary choice of the two complex
numbers $a$, $b$ in (A.1c), and for an appropriate choice of the {\it
real} quantity $\Omega t$, the solution (A.1c) of (A.1a) is {\it
real}. Indeed (A.1c) entails
\begin{gather}
{\rm Re} \,(x) = [{\rm Re}\, (a) + {\rm Re}\, (b) ]
\cos (\Omega t) + [{\rm Im}\,
(a) - {\rm Im}\, (b)]\sin (\Omega t),\\
{\rm Im}\, (x) =[ -{\rm Re}\,(a) + {\rm Re}\,(b) ] \sin (\Omega t) +
[{\rm Im}\, (a) + {\rm Im} \,(b)] \cos
(\Omega t).
\end{gather}
\end{subequations}
Now choose $t$ (${\rm mod}\, (T)$, see (1.47)) so that
\begin{subequations}
\begin{gather}
\tan (\Omega t) = [{\rm Im}\, (a) + {\rm Im}\, (b)] /[{\rm Re}\, (a)
-{\rm Re}\,(b)],
\end{gather}
entailing
\begin{gather}
\sin (\Omega t) = [{\rm Im} \,(a) + {\rm Im}\, (b)]/D,\\
\cos (\Omega t) = [{\rm Re} \, (a) - {\rm Re} \, (b)]/D,\\
D=\left\{ [{\rm Im}\, (a) + {\rm Im}\, (b)]^{2} +
[{\rm Re}\, (a) - {\rm Re}\, (b)]^{2}\right\}^{1/2}.
\end{gather}
\end{subequations}
Note that a {\it real} choice of $t$ such that (A.2a) hold can {\it
always} be made; and, via (A.1f), this choice entails that ${\rm Im}\, (x)$
vanishes, namely that $x$ is {\it real}; indeed $x$ is then given, via
(A.2b), (A.2c), (A.2d), by the explicit expression
\begin{gather}
x = \left(\vert a \vert^{2} - \vert b \vert^{2}\right) / D.
\end{gather}
This completes our proof.

Next we prove that the inequality (2.57) is (of course for $\Omega$
{\it real}) sufficient to guarantee that the solution $w (x,t)$, see
(2.56), is {\it nonsingular} for all {\it real} values of $x$ and $t$,
namely that the two roots
\begin{gather}
x_{s} (t) = a \exp (- {\rm i} \Omega t) + s [{\rm i} (\beta / \Omega)
+ b \exp (-2{\rm i} \Omega t)]^{1/2}, \qquad
s= \pm,
\end{gather}
of the polynomial in $x$ in the denominator in the right-hand side of
(2.56) have nonvanishing imaginary parts for {\it all} ({\it real})
values of the time $t$, if the inequality (2.57) holds.

To prove this, we note first of all that the reality of $\Omega$ and
$t$ obviously entails
\begin{subequations}
\begin{gather}
\left\vert {\rm Im}\, [{\rm i}(\beta / \Omega) +
b \exp (-2 {\rm i} \Omega t)] \right\vert \geq
\vert {\rm Re} \, (\beta) / \Omega \vert - \vert b \vert,
\end{gather}
hence (if $\vert {\rm Re}\, (\beta)/ \Omega \vert > \vert b \vert$,
as indeed implied by (2.57))
\begin{gather}
\left\vert {\rm Im}\, [{\rm i} (\beta /\Omega) + b \exp (-2 {\rm i}
\Omega t)]^{1/2} \right\vert >
(\vert {\rm Re}\, (\beta ) /\Omega \vert - \vert b
\vert)^{1/2}/2,
\end{gather}
hence
\begin{gather}
\left\vert {\rm Im}\,\left\{[a \exp (- {\rm i} \Omega t)] \pm [{\rm
i}(\beta /\Omega) + b \exp (-2 {\rm i} \Omega t)]^{1/2}\right\}
\right\vert\nonumber\\
\qquad {}> (\vert {\rm Re}\, (\beta) / \Omega \vert - \vert
b \vert)^{1/2} /2 - \vert a \vert.
\end{gather}

This entails, via (A.4),
\begin{gather}
\vert {\rm Im} \, [x_{s} (t)]\vert > (\vert {\rm Re} \, (\beta)
/\Omega \vert - \vert b \vert)^{1/2} /2 - \vert a \vert,\qquad s=
\pm,
\end{gather}
\end{subequations}
hence, if (2.57) holds,
\begin{gather}
\vert {\rm Im} \, [x_{s} (t) ]\vert >0,\qquad s = \pm,
\end{gather}
QED.

\label{calogero-lastpage}
\end{document}